\documentclass[12pt]{amsart}  	

\usepackage[
  left=0.85in,
  right=0.85in,
  top=0.9in,
  bottom=0.9in
]{geometry}  
\usepackage{comment}
\usepackage[all]{xy}						
\usepackage{amssymb}
\usepackage{graphicx}
\usepackage{amscd,latexsym,amsthm,amsfonts,amssymb,amsmath,amsxtra}
\usepackage[mathscr]{eucal}
\usepackage{relsize} 
\usepackage[colorlinks]{hyperref}
\usepackage{mathtools}
\usepackage{multirow}
\usepackage{multicol}
\usepackage{empheq}
\usepackage{enumitem}
\usepackage{commath}
\usepackage{subfiles}
\usepackage{quiver}
\usepackage{chngcntr}
\numberwithin{equation}{section}
\usepackage{bbm}
\newcounter{keepeqno}

\usepackage{tikz}
\usetikzlibrary{decorations.pathreplacing}
\usepackage{adjustbox}

\hypersetup{linkcolor=black, citecolor=black}

\newcommand{\cA}{{\mathcal{A}}}   
\newcommand{\cC}{{\mathcal{C}}}   \newcommand{\cD}{{\mathcal{D}}}
   \newcommand{\cF}{{\mathcal{F}}}
   \newcommand{\cH}{{\mathcal{H}}}

\newcommand{\cO}{{\mathcal{O}}}

\newcommand{\cW}{{\mathcal{W}}}

\newcommand{\St}{{\mathrm {St}}}

\newcommand\restr[2]{{
  \left.\kern-\nulldelimiterspace 
  #1 
  \vphantom{\big|} 
  \right|_{#2} 
  }}

\newcommand{\supp}{{\mathrm{Supp}}}

\newcommand{\Irr}{{\mathrm{Irr}}}

\newcommand{\GL}{{\mathrm{GL}}}

\newcommand{\Hom}{{\mathrm{Hom}}}

\newcommand{\Ind}{{\mathrm{Ind}}}
\newcommand{\cInd}{{\mathrm{c-Ind}}}
\newcommand{\uInd}{{\mathrm{u-Ind}}}

\newcommand{\Jac}{{\mathrm{Jac}}}

\newcommand{\Mat}{{\mathrm{Mat}}}

\renewcommand{\Re}{{\mathrm{Re}}}

\newcommand{\Sl}{{\mathrm{SL}}}

\newcommand{\Sp}{{\mathrm{Sp}}}
\newcommand{\Stab}{{\mathrm{Stab}}}

\newcommand{\Span}{{\mathrm{Span}}}
\newcommand{\one}{\mathbf{1}}

\newcommand{\w}{\text{\textbf{w}}}

\def\diag{{\rm diag}}

\def\Rep{\mathrm{Rep}}

\theoremstyle{plain}
\newtheorem{thm}{Theorem}[section]
\newcommand{\restatedtheoremheading}{}
\newtheorem*{restatedtheoreminner}{\restatedtheoremheading}
\newenvironment{restatedthm}[1]
  {\renewcommand{\restatedtheoremheading}{Theorem~\ref{#1}}\begin{restatedtheoreminner}}
  {\end{restatedtheoreminner}}
\newtheorem{prp}[thm]{Proposition}
\newtheorem{lem}[thm]{Lemma}
\newtheorem{cor}[thm]{Corollary}
\newtheorem{cnj}[thm]{Conjecture}
\newtheorem*{restatedconjectureinner}{Conjecture~\ref{cnj:odd-tempered-necessity}}

\theoremstyle{definition}

\newtheorem{exm}[thm]{Example}

\theoremstyle{remark}
\newtheorem{rmk}[thm]{\textbf{Remark}}

\makeatletter

\newcommand{\Rmnum}[1]{\expandafter\@slowromancap\romannumeral #1@}
\makeatother

\DeclareMathOperator{\Gl}{\mathrm{GL}}

\DeclarePairedDelimiter{\floor}{\lfloor}{\rfloor}

\title{Plancherel Decompositions and Distinctions for $\Sp_{2n}\backslash \Gl_{2n+1}$ and $P(\Sp_{2n})\backslash \Gl_{2n}$}

\author{Sagnik Mukherjee}
\address{School of Mathematics, University of Minnesota,
206 Church St. S.E., Minneapolis, MN 55455, USA}

\makeatletter
\def\@author{\authors}
\makeatother

\begin{document}
\begin{abstract}
Let $F$ be a non-Archimedean local field of characteristic zero.
We determine the Plancherel decompositions of
\(
\Sp_{2n}(F)\backslash\GL_{2n+1}(F)
\ \text{and}\
P(\Sp_{2n})(F)\backslash\GL_{2n}(F),
\)
where $P(\Sp_{2n})$ denotes the mirabolic subgroup of $\Sp_{2n}(F)$, and investigate smooth distinction for both
spaces. We carry out the local unfolding procedure of
\cite[\S9.5]{SakellaridisVenkatesh} explicitly for these spaces, constructing
unitary isomorphisms with degenerate Whittaker models by successive Fourier
transforms. The isomorphism for the even model was already stated in
\cite[Example~9.5.17]{SakellaridisVenkatesh}. The resulting Plancherel
formulas are expressed as direct integrals of normalized
parabolic inductions from the tempered duals of
$\GL_n(F)\times\GL_{n+1}(F)$ and
$\GL_n(F)\times\GL_n(F)$, respectively, with the corresponding
product Whittaker--Plancherel measures. In the odd case,
the parameter map agrees with the prediction of
\cite{BenZviSakellaridisVenkatesh}.

Using closed parabolic orbits and local
Rankin--Selberg integrals, we prove that every representation
induced from these tempered data admits a nonzero smooth
invariant functional for the corresponding model.
We also prove that essentially square-integrable
representations of $\GL_{2n+1}(F)$ are not
$\Sp_{2n}(F)$-distinguished, and that supercuspidal
representations of $\GL_{2n}(F)$ are not
$P(\Sp_{2n})(F)$-distinguished.
For the even model with $n=1,2$, we construct
$P(\Sp_{2n})(F)$-distinguished square-integrable
representations of $\GL_{2n}(F)$ lying outside the
Plancherel support. Thus the smooth tempered distinguished
spectrum can be strictly larger than the corresponding
$L^2$-spectrum. For $\GL_5(F)$, however, we prove that an irreducible
tempered representation is $\Sp_4(F)$-distinguished if and only if it is
induced from tempered representations of $\GL_2(F)$ and $\GL_3(F)$.
\end{abstract}

\maketitle 
\section{Introduction}

Let \(F\) be a non-Archimedean local field of characteristic zero. Given a reductive group \(G\) over \(F\) and a closed unimodular subgroup \(H\subset G(F)\), two basic questions arise in the study of the homogeneous space \(H\backslash G(F)\). The first asks which irreducible smooth representations \(\pi\) of \(G(F)\) admit a nonzero \(H\)-invariant linear form, that is,
\[
\Hom_H(\pi,\mathbf 1)\neq 0.
\]
Such representations are called \(H\)-distinguished. The second asks for the Plancherel decomposition of \(L^2(H\backslash G(F))\), including its spectral measure and multiplicities. Writing \(X=H\backslash G(F)\), this decomposition under $G$-right translation takes the form \cite[\S\S3.1.2--3.1.3, (3.6)]{BeuzartPlessisRelativeLanglands}\[L^2(X)\simeq\int_{\Pi_{\mathrm{unit}}(G)}^\oplus\mathcal M_X(\pi)\widehat\otimes\mathcal H_\pi\,d\mu_X(\pi).\]Here \(\Pi_{\mathrm{unit}}(G)\) is the unitary dual of \(G(F)\), \(\mathcal H_\pi\) is the Hilbert space of \(\pi\), and \(\mathcal M_X(\pi)\) is a Hilbert multiplicity space on which \(G(F)\) acts trivially; \(\widehat\otimes\) denotes the Hilbert tensor product. The measure \(\mu_X\) represents the Plancherel measure class, and \(\dim\mathcal M_X(\pi)\) is the spectral multiplicity for \(\mu_X\)-almost every \(\pi\). These questions are closely related to the study of period integrals of automorphic forms. A central aim of the relative Langlands program is to describe the representations that arise in these problems in terms of their Langlands or Arthur parameters \cite{SakellaridisVenkatesh,BenZviSakellaridisVenkatesh}.

In this paper, we study the Plancherel decomposition and smooth distinction for the two spaces
\[
\begin{aligned}
X_n^{\mathrm{odd}}
   =\Sp_{2n}(F)\backslash\GL_{2n+1}(F),\qquad
X_n^{\mathrm{even}}
   =P(\Sp_{2n})\backslash\GL_{2n}(F),
\end{aligned}
\]
where \(n\geq 1\). Here \(\Sp_{2n}(F)\) is embedded in \(\GL_{2n+1}(F)\) by \(h\mapsto\operatorname{diag}(h,1)\), and
\[
P(\Sp_{2n})
   =\{h\in\Sp_{2n}(F):he_1=e_1\}
\]
is its Klingen mirabolic subgroup, where \(e_1\) is the first standard basis vector of \(F^{2n}\). 

Our study is motivated in part by the classical symplectic model
\[
\Sp_{2n}(F)\backslash\GL_{2n}(F).
\]
The local theory began with models over finite fields \cite{KlyachkoModels} and their analogues over non-Archimedean local fields \cite{HR90}. In particular, an irreducible generic representation of \(\GL_{2n}(F)\) cannot admit a symplectic period \cite[Theorem~3.2.2]{HR90}. The symplectic and Whittaker models are the two extremes of the family of Klyachko models. Every irreducible unitary representation admits one of these models \cite{OffSayMixed}, and the models satisfy uniqueness and disjointness \cite{OffSayDisj}. Their behavior under parabolic induction and the classification for ladder representations are studied in \cite{OffMitSay}. The global theory concerns symplectic and mixed periods on the discrete automorphic spectrum of general linear groups \cite{JacquetRallisSymplectic,OffenSymplecticPeriods,OffSayMixed}. For an account of the relation between period integrals and local distinction, see \cite{Off}.

The pair $(\GL_r(D),\Sp_r(D))$, with \(D/F\) the quaternion division algebra, is considered in \cite{VermaTemperedSymplectic}. The distinguished irreducible tempered representations are classified in \cite[Theorem~1.1]{VermaTemperedSymplectic}, proving the revised conjecture recorded in \cite[Conjecture~4.1]{VermaTemperedSymplectic}. They are precisely the normalized parabolic inductions of factors \(\rho\) or \(\delta_{m+1}(\rho)\times\delta_m(\rho)\), where \(m\geq1\) and \(\rho\) is unitary, supercuspidal, and symplectically distinguished. Here \(\delta_m(\rho)\) denotes the generalized Steinberg representation of length \(m\) attached to \(\rho\).

The Plancherel formula in \cite{LapidOffenPlancherel} expresses
\[
L^2\bigl(\Sp_{2n}(F)\backslash\GL_{2n}(F)\bigr)
\]
as a multiplicity-free direct integral of Speh representations \(U(\tau,2)\), with \(\tau\) ranging over the tempered dual of \(\GL_n(F)\), equipped with its group Plancherel measure. Here \(U(\tau,2)\) is the unique irreducible quotient of
\[
\nu^{1/2}\tau\times\nu^{-1/2}\tau,
\]
where \(\nu(g)=|\det g|_F\) and \(\times\) denotes normalized parabolic induction. This formula provides a natural point of comparison for our two spaces, obtained by enlarging the ambient general linear group by one dimension or by replacing the symplectic group with its vector stabilizer.

A further motivation comes from relative Langlands duality. For spherical homogeneous spaces, the conjectures in \cite{SakellaridisVenkatesh} describe the Plancherel spectrum through distinguished Arthur parameters associated with a dual group attached to the space. This picture is extended to a duality between Hamiltonian varieties in \cite{BenZviSakellaridisVenkatesh}. In the case of a cotangent bundle \(T^*X\), this relates the harmonic analysis of \(X\) to spectral data attached to the dual Hamiltonian variety.

For
\(
\mathfrak M_n
   =T^*\bigl(\Sp_{2n}\backslash\GL_{2n+1}\bigr)
\)
(the polarised case), the dual predicted in \cite{BenZviSakellaridisVenkatesh} is
\(
\check{\mathfrak M}_n
   =T^*\bigl((\GL_n\times\GL_{n+1})\backslash\GL_{2n+1}\bigr).
\)
Thus the symplectic period on \(\GL_{2n+1}\) is paired with the linear model associated with \(\GL_n\times\GL_{n+1}\). This suggests that the local Plancherel decomposition of \(X_n^{\mathrm{odd}}\) should involve the block-diagonal parameter map
\[
\GL_n(\mathbb C)\times\GL_{n+1}(\mathbb C)
   \longrightarrow \GL_{2n+1}(\mathbb C).
\]
This dual pair is studied globally in \cite{Xi-Lu}. The vanishing of cuspidal symplectic periods \cite{Ash-Ginzburg-Rallis} is extended there to Eisenstein series induced from cuspidal data on maximal Levi subgroups: at regular parameters the period vanishes unless the Levi is associate to \(\GL_n\times\GL_{n+1}\) \cite[Theorem~6.3.1]{Xi-Lu}. For this Levi, an Eulerian period formula is obtained in \cite[Theorem~1.1.3]{Xi-Lu}. The corresponding cuspidal vanishing in the geometric setting is proved in \cite[Appendix~A]{Xi-Lu}.

Motivated by these results, we study the local Plancherel decompositions and smooth distinction for \(X_n^{\mathrm{odd}}\) and \(X_n^{\mathrm{even}}\). We are particularly interested in how the representations appearing in these decompositions compare with those admitting invariant linear forms. Our \(L^2\)-calculation follows the local unfolding procedure of \cite[\S9.5, Theorem~9.5.9]{SakellaridisVenkatesh}, which applies Fourier transforms along abelian unipotent quotients. The isomorphism for \(X_n^{\mathrm{even}}\) is already stated in \cite[Example~9.5.17]{SakellaridisVenkatesh}. Here we give an explicit verification of successive unfolding for the two models, determining the intermediate subgroups and characters and constructing the unitary maps from isometries on dense subspaces. This also yields the isomorphism for \(X_n^{\mathrm{odd}}\) with a degenerate Whittaker model. The resulting Plancherel decompositions then follow from the Whittaker--Plancherel formula \cite{DelormeWhittakerPlancherel}. Finally, we prove by a purely local argument that no irreducible essentially square-integrable representation of \(\mathrm{GL}_{2n+1}(F)\) is \(\mathrm{Sp}_{2n}(F)\)-distinguished. However, for the model \((\mathrm{GL}_{2n}(F),P(\mathrm{Sp}_{2n}(F)))\), the analogous statement fails: we construct irreducible square-integrable representations with nonzero invariant linear forms for \(n=1\) and \(n=2\). We also determine exactly which irreducible tempered representations of \(\GL_5(F)\) are \(\Sp_4(F)\)-distinguished. We state these results below.

\subsection{Main results}
Fix a nontrivial unitary additive character \(\psi:F\to\mathbb C^\times\). Let \(N_r\) be the standard maximal upper triangular unipotent subgroup of \(\GL_r(F)\), and put
\[
\psi_r(u)=\psi(u_{1,2}+\cdots+u_{r-1,r}),
\qquad \psi_1=1.
\]
For \(a+b=r\), define
\[
\psi_{a,b}
\begin{pmatrix}
u_1&*\\
0&u_2
\end{pmatrix}
   =\psi_a(u_1)\psi_b(u_2),
\qquad u_1\in N_a,\quad u_2\in N_b.
\]
This character is trivial on the unipotent radical of the standard parabolic \(P_{a,b}\). We write \(L^2(N_r\backslash\GL_r(F),\psi_{a,b})\) for the Hilbert space of functions with this left equivariance that are square-integrable modulo \(N_r\), with the right translation action.

Our first result identifies the two \(L^2\)-spaces with degenerate Whittaker models.

\begin{thm}\label{thm:intro-unitary-unfolding}
There are unitary equivariant isomorphisms
\[
L^2\bigl(\Sp_{2n}(F)\backslash\GL_{2n+1}(F)\bigr)
   \simeq
L^2\bigl(N_{2n+1}\backslash\GL_{2n+1}(F),\psi_{n,n+1}\bigr)
\]
and
\[
L^2\bigl(P(\Sp_{2n})\backslash\GL_{2n}(F)\bigr)
   \simeq
L^2\bigl(N_{2n}\backslash\GL_{2n}(F),\psi_{n,n}\bigr).
\]
\end{thm}

Unitary induction in stages \cite[Chapter~III, Theorem~3.3]{Mackey} and the Whittaker--Plancherel formula \cite[Theorem~7.7]{DelormeWhittakerPlancherel} give the following decompositions. Here \(\Pi_{\mathrm{temp}}(\GL_r)\) is the tempered unitary dual of \(\GL_r(F)\), \(d\mu_r\) is its Whittaker--Plancherel measure, and \(\mathcal H_\pi\) is the Hilbert space of \(\pi\).

\begin{thm}\label{thm:intro-plancherel}
As unitary representations,
\[
\begin{aligned}
&L^2\bigl(\Sp_{2n}(F)\backslash\GL_{2n+1}(F)\bigr)\\
&\qquad\simeq
\int_{\Pi_{\mathrm{temp}}(\GL_n)
      \times\Pi_{\mathrm{temp}}(\GL_{n+1})}^{\oplus}
\mathcal H_{\pi_n\times\pi_{n+1}}\,
d\mu_n(\pi_n)\,d\mu_{n+1}(\pi_{n+1}),
\end{aligned}
\]
and
\[
\begin{aligned}
&L^2\bigl(P(\Sp_{2n})\backslash\GL_{2n}(F)\bigr)\\
&\qquad\simeq
\int_{\Pi_{\mathrm{temp}}(\GL_n)
      \times\Pi_{\mathrm{temp}}(\GL_n)}^{\oplus}
\mathcal H_{\pi_1\times\pi_2}\,
d\mu_n(\pi_1)\,d\mu_n(\pi_2).
\end{aligned}
\]
\end{thm}

All the induced representations in these formulas are irreducible and tempered \cite{Bern-ZeleI,ZelII}. Their Langlands parameters are, respectively,
\[
\phi_{\pi_n}\oplus\phi_{\pi_{n+1}},
\qquad
\phi_{\pi_1}\oplus\phi_{\pi_2}.
\]
In the odd case, the parameter map therefore agrees with the prediction described above. The formulas retain the inducing data: distinct pairs can give isomorphic representations of the ambient general linear group. The spectral multiplicities are obtained by disintegrating the product measures along the induction maps.

We next consider smooth distinction. The existence of an invariant functional realizes a representation in \(C^\infty(H\backslash G)\), but does not by itself place it in the Plancherel support of \(L^2(H\backslash G)\). Conversely, the Plancherel decomposition alone does not establish smooth distinction at every individual parameter. We prove the following statement for all the tempered inducing data in Theorem~\ref{thm:intro-plancherel}.

\begin{thm}\label{thm:intro-sufficient}
Let \(n\geq1\).
\begin{enumerate}[label=\textup{(\arabic*)},leftmargin=2em]
\item If \(\sigma_1\) and \(\sigma_2\) are irreducible tempered representations of \(\GL_n(F)\) and \(\GL_{n+1}(F)\), respectively, then \(\sigma_1\times\sigma_2\) is \(\Sp_{2n}(F)\)-distinguished.
\item If \(\tau_1\) and \(\tau_2\) are irreducible tempered representations of \(\GL_n(F)\), then \(\tau_1\times\tau_2\) is \(P(\Sp_{2n})\)-distinguished.
\end{enumerate}
\end{thm}

For the odd model, we prove that no irreducible essentially square-integrable representation is distinguished.

\begin{thm}\label{thm:intro-discrete-vanishing}
If \(\pi\) is an irreducible essentially square-integrable representation of \(\GL_{2n+1}(F)\), then
\[
\Hom_{\Sp_{2n}(F)}(\pi,\mathbf1)=0.
\]
\end{thm}

For \(P(\Sp_{2n})\), the argument in \cite[Proposition~3, Corollary~1 and Remark~2]{DijolsPrasad} gives
\[
\Hom_{P(\Sp_{2n})}(\pi,\mathbf1)=0
\]
for every irreducible supercuspidal representation \(\pi\) of \(\GL_{2n}(F)\). We recall this argument in our coordinates in Proposition~\ref{prp:even-supercuspidal-vanishing}.

This vanishing for \(P(\Sp_{2n})\) does not extend to square-integrable representations. For \(n=1\), one has \(P(\Sp_2)=N_2\), and every unitary twist of the Steinberg representation of \(\GL_2(F)\) is \(N_2\)-distinguished. For \(\GL_4(F)\), we obtain the following family.

\begin{thm}\label{thm:intro-even-counterexample}
Let \(\rho\) be an irreducible unitary supercuspidal representation of \(\GL_2(F)\), and put
\[
\delta_\rho
   =\operatorname{St}
      \bigl([\nu^{-1/2}\rho,\nu^{1/2}\rho]\bigr),
\]
where \(\operatorname{St}(\Delta)\) denotes the essentially square-integrable representation attached to the Zelevinsky segment \(\Delta\) \cite{ZelII}. Then \(\delta_\rho\) is an irreducible square-integrable representation of \(\GL_4(F)\), and
\[
\dim\Hom_{P(\Sp_4)}(\delta_\rho,\mathbf1)=1.
\]
\end{thm}

These square-integrable representations of \(\GL_2(F)\) and \(\GL_4(F)\) lie outside the respective Plancherel supports in Theorem~\ref{thm:intro-plancherel}. Thus, even among tempered representations, smooth distinction for the even model is not exhausted by the inducing data in its Plancherel formula.

The counterexample in Theorem~\ref{thm:intro-even-counterexample} raises the question whether the sufficient condition in Theorem~\ref{thm:intro-sufficient} is also necessary for \((\GL_5(F),\Sp_4(F))\). To see the issue, let \(\chi\) be a unitary character of \(F^\times\). Then \(\chi\times\delta_\rho\) is an irreducible tempered representation of \(\GL_5(F)\) that is not induced from tempered data on \(\GL_2(F)\times\GL_3(F)\). Let \(V_{\mathrm{op}}\) be the subrepresentation of \((\chi\times\delta_\rho)|_{\Sp_4(F)}\) consisting of sections compactly supported on the open orbit in \(P_{1,4}(F)\backslash\GL_5(F)\). This orbit has stabilizer \(P(\Sp_4)\), so the geometric lemma (Lemma~\ref{lem: geom lem II}) gives
\[
 \Hom_{\Sp_4(F)}(V_{\mathrm{op}},\mathbf1)
 \simeq \Hom_{P(\Sp_4)}(\delta_\rho,\mathbf1).
\]
The functional in Theorem~\ref{thm:intro-even-counterexample} therefore gives a nonzero invariant functional on the open-orbit subspace. An extension to \(\chi\times\delta_\rho\) would contradict the proposed necessary condition for the odd model. We prove that no such extension exists. More generally, we establish the necessary condition for every irreducible tempered representation of \(\GL_5(F)\). Combined with Theorem~\ref{thm:intro-sufficient}, this gives the following characterization.

\begin{thm}\label{thm:intro-GL5-necessary}
Let \(\pi\) be an irreducible tempered representation of \(\GL_5(F)\). Then
\[
 \Hom_{\Sp_4(F)}(\pi,\mathbf1)\neq0
 \quad\Longleftrightarrow\quad
 \pi\simeq\sigma_2\times\sigma_3
\]
for some irreducible tempered representations \(\sigma_2\) of
\(\GL_2(F)\) and \(\sigma_3\) of \(\GL_3(F)\).
\end{thm}

Thus, Theorem~\ref{thm:intro-GL5-necessary} shows that the additional distinguished representations for the even model do not give counterexamples for the odd model when \(n=2\). For \(n=1\), the same necessary condition follows from Theorem~\ref{thm:intro-discrete-vanishing} and the tempered classification: every nonsquare-integrable tempered representation of \(\GL_3(F)\) can be written as \(\sigma_1\times\sigma_2\). These cases, together with Theorem~\ref{thm:intro-sufficient}, suggest the following conjecture.

\begin{cnj}\label{cnj:odd-tempered-necessity}
Let \(n\geq1\), and let \(\pi\) be an irreducible tempered representation of \(\GL_{2n+1}(F)\). If \(\Hom_{\Sp_{2n}(F)}(\pi,\mathbf1)\neq0\), then
\[
 \pi\simeq\sigma_n\times\sigma_{n+1}
\]
for irreducible tempered representations \(\sigma_n\) of \(\GL_n(F)\) and \(\sigma_{n+1}\) of \(\GL_{n+1}(F)\).
\end{cnj}

We briefly describe the proofs. The proof of Theorem~\ref{thm:intro-unitary-unfolding} is a local unfolding: a sequence of Fourier transforms along abelian unipotent quotients. At each stage, the nonzero characters form a single symplectic orbit whose complement has measure zero. Fourier--Plancherel \cite{PoonenTate} therefore replaces the original model by one attached to the stabilizer of a nonzero character. The first such stabilizer in the odd model is \(P(\Sp_{2n})\), explaining the appearance of the even model. Repeating the construction in smaller symplectic spaces, and applying unitary induction in stages \cite[Chapter~III, Theorem~3.3]{Mackey}, gives Theorem~\ref{thm:intro-unitary-unfolding}.

The smooth invariant functionals in Theorem~\ref{thm:intro-sufficient} are constructed from the closed identity orbits for \(P_{n,n+1}\) in the odd model and \(P_{n,n}\) in the even model. The required forms on the inducing representations come from local Rankin--Selberg periods \cite{MR701565,Jo-Local-Periods}; the closed-orbit argument of \cite[Lemma~4.3]{OffenParabolicInduction}, in the form proved in Proposition~\ref{prp: closed orbit contribution}, gives functionals on the full induced representations. The proof of Theorem~\ref{thm:intro-discrete-vanishing} combines the Kirillov model \cite{Bern-ZeleI} and the distribution vanishing statement in \cite[Proposition~1]{OffSayDisj} with the restriction theorem in \cite[Theorem~5.7]{Chan-Savin}. The counterexample for \(\GL_4(F)\) follows from smooth Clifford theory \cite[Proposition~1]{DijolsPrasad} and the Jacquet modules of \(\delta_\rho\). For Theorem~\ref{thm:intro-GL5-necessary}, the geometric lemma reduces the remaining distinction problem to invariant functionals on an affine chart. We show that such a functional would satisfy incompatible scaling relations under dilations and Fourier transform.

\subsection{Organization of the paper}
\label{subsec:intro-organization}

Section~2 fixes notation and recalls the representation-theoretic results used below. Section~\ref{sec:local-unfolding} proves the \(L^2\)-isomorphisms and the Plancherel formulas. Section~\ref{sec:sufficient-distinction} proves the sufficient conditions for distinction in Theorem~\ref{thm:intro-sufficient}, using closed parabolic orbits and local Rankin--Selberg periods. Section~\ref{sec:limits-smooth-necessity} proves the vanishing results of Theorem~\ref{thm:intro-discrete-vanishing} and Proposition~\ref{prp:even-supercuspidal-vanishing}, and gives square-integrable examples for \(\GL_2(F)\) and \(\GL_4(F)\) showing that smooth distinction can occur outside the Plancherel support. It concludes with the characterization of tempered \(\Sp_4(F)\)-distinguished representations of \(\GL_5(F)\) in Theorem~\ref{thm:intro-GL5-necessary}. Appendix~\ref{app:parabolic-orbits} records the orbits, stabilizers, and Levi projections for
\[
P_{\ell,\,2n+1-\ell}\backslash
\GL_{2n+1}/\Sp_{2n},
\qquad 1\leq\ell\leq n,
\]
extending the calculations for \(P_{n,n+1}\) in \cite[\S5]{Xi-Lu} to these maximal parabolics.

\section*{Acknowledgements and AI-Disclosure}
I am grateful to Chen Wan for suggesting this problem and for guiding me in obtaining the Plancherel decomposition by local unfolding. I thank my advisor, Dihua Jiang, for his guidance and many useful discussions, and especially for suggesting the use of Bernstein--Zelevinsky derivatives to study smooth distinction. I also thank Sarbartha Bhattacharya, Zhaolin Li, David Luo, and Guodong Xi for many helpful conversations.

During the final stages of this project, the author used ChatGPT (Sol and Astra) to discuss the arguments, assist in checking mathematical claims, improve the exposition, help with \LaTeX, and locate references. Nevertheless, the author takes full responsibility for the mathematical results, citations, and final text of this paper.
\par\vspace{0.5cm}
\section{Notation and preliminaries}
\subsection{Notation and conventions}
\label{subsec:notation}

Let \(F\) be a non-Archimedean local field of characteristic zero. All
representations are complex and smooth unless explicitly realized on a
Hilbert space. In unitary direct integrals we use their Hilbert
completions; spaces \(\Hom_H\) refer to equivariant linear maps on the
underlying smooth representations. For a locally compact totally
disconnected group \(L\), write \(\Rep(L)\) for its category of smooth
representations, \(\Irr(L)\) for the set of irreducible isomorphism
classes, and \(\one_L\) for the trivial representation. The smooth
contragredient of \(\pi\) is denoted by \(\widetilde\pi\). For a closed
subgroup \(L'\subseteq L\) and a character \(\chi\) of \(L'\), we say
that \(\pi\) is \((L',\chi)\)-distinguished if
\(\Hom_{L'}(\pi,\chi)\neq0\), and simply \(L'\)-distinguished when
\(\chi=\one_{L'}\).

We often suppress \(F\) from the notation for matrix groups and put
\(\Gl_0=\{1\}\). Let \(I_r\) be the identity matrix,
\(\Mat_{a,b}(F)\) the space of \(a\times b\) matrices, and
\(e_1,\ldots,e_r\) the standard basis of \(F^r\). Write \(\w_{0,r}\)
for the matrix with ones on the anti-diagonal and zeros elsewhere.
For every \(r\), set \(\nu(g)=\abs{\det g}_F\) on \(\Gl_r\) and
\(\nu^s\pi=\pi\otimes\nu^s\).

For a composition \(\mathbf n=(n_1,\ldots,n_t)\) of \(r\), let
\(P_{\mathbf n}=M_{\mathbf n}\ltimes U_{\mathbf n}\) be the standard
upper triangular parabolic subgroup of \(\Gl_r\), with
\(M_{\mathbf n}\simeq\Gl_{n_1}\times\cdots\times\Gl_{n_t}\).
The modular character of a locally compact group \(L\) is denoted by
\(\delta_L\). If \(P=MU\) is a parabolic subgroup of a reductive group
\(G\), we inflate \(\sigma\in\Rep(M)\) to \(P\) and use normalized
parabolic induction:
\[
 \Ind_P^G(\sigma)
 =\{f:G\to V_\sigma\text{ smooth}:
 f(mug)=\delta_P(m)^{1/2}\sigma(m)f(g)\},
\]
with \(G\) acting by right translation. For representations
\(\sigma_i\) of \(\Gl_{n_i}\), write
\[
 \sigma_1\times\cdots\times\sigma_t
 :=\Ind_{P_{\mathbf n}}^{\Gl_r}
 (\sigma_1\boxtimes\cdots\boxtimes\sigma_t).
\]
We reserve \(\cInd\) for non-normalized compact induction.

For a unipotent subgroup \(U\) and a character \(\chi\) of \(U\),
write \(\pi_{U,\chi}\) for the unnormalized twisted coinvariant space
\[
 V_\pi/\Span\{\pi(u)v-\chi(u)v:u\in U,\ v\in V_\pi\},
\]
and put \(\pi_U=\pi_{U,\one}\). The stabilizer of \(\chi\) in the
normalizer of \(U\) acts on this quotient through \(\pi\). For a
parabolic subgroup \(P=MU\), the normalized Jacquet module is
\(\Jac_P^G(\pi)=\delta_P^{-1/2}\otimes\pi_U\), regarded as an
\(M\)-module. We also write \(\Jac_{(n_1,\ldots,n_t)}(\pi)\) when
\(G=\Gl_r\) and \(P=P_{n_1,\ldots,n_t}\).

We realize the symplectic group in the ordered basis
\(e_1,\ldots,e_r,f_r,\ldots,f_1\) by
\[
 J_{2r}=\begin{pmatrix}0&\w_{0,r}\\-\w_{0,r}&0\end{pmatrix},
 \qquad
 \Sp_{2r}=\{h\in\Gl_{2r}:{}^thJ_{2r}h=J_{2r}\}.
\]
Its mirabolic subgroup is \(P(\Sp_{2r})=\Stab_{\Sp_{2r}}(e_1)\).
For \(g\in\Gl_r\), put
\(g^\#=\w_{0,r}\,{}^tg^{-1}\w_{0,r}\). This is an involutive
automorphism, and we write \(\sigma^\#(g)=\sigma(g^\#)\).

\subsection{Whittaker and Kirillov models}

Fix a nontrivial additive character \(\psi:F\to\mathbb C^\times\).
Let \(N_r\) be the standard maximal upper triangular unipotent
subgroup of \(\Gl_r\), equipped with the character
\[
 \psi_r(u)=\psi\!\left(\sum_{i=1}^{r-1}u_{i,i+1}\right).
\]
An irreducible representation \(\pi\) of \(\Gl_r\) is generic if
\(\Hom_{N_r}(\pi,\psi_r)\neq0\); this space is then one-dimensional
\cite[Theorem~9.7]{ZelII}. Choosing a nonzero Whittaker functional
\(\ell\) realizes \(\pi\) in its Whittaker model
\(\mathcal W(\pi,\psi_r)=\{g\mapsto\ell(\pi(g)v):v\in V_\pi\}\).

Let \(Q_1=\{1\}\), and for \(r\geq2\) let
\[
 Q_r=\left\{\begin{pmatrix}g&x\\0&1\end{pmatrix}:
 g\in\Gl_{r-1},\ x\in F^{r-1}\right\}
\]
be the mirabolic subgroup of \(\Gl_r\). We need the following standard
form of the Kirillov model, which follows from
\cite[\S3.5 and Theorem~4.4]{Bern-ZeleI}: the mirabolic filtration has
only its highest-derivative term, since all intermediate derivatives of a
supercuspidal representation vanish. Its highest derivative is
one-dimensional.

\begin{prp}\label{prp:supercuspidal-restriction-mirabolic}
If \(\rho\) is an irreducible supercuspidal representation of
\(\Gl_r\), with \(r\geq2\), then
\[
 \restr{\rho}{Q_r}\simeq\cInd_{N_r}^{Q_r}(\psi_r),
 \qquad
 \restr{\rho}{\Gl_{r-1}}
 \simeq\cInd_{N_{r-1}}^{\Gl_{r-1}}(\psi_{r-1}),
\]
where \(\Gl_{r-1}\) is embedded in the upper-left block.
\end{prp}

The second isomorphism is obtained from the first by restricting
functions to \(\Gl_{r-1}\), since \(Q_r=N_r\Gl_{r-1}\) and
\(N_r\cap\Gl_{r-1}=N_{r-1}\).

\subsection{Segments and Jacquet modules}

Let \(\rho\in\Irr(\Gl_d)\) be supercuspidal, and let \(a,b\in\mathbb C\)
satisfy \(b-a\in\mathbb Z_{\geq0}\). For the segment
\(\Delta=[\nu^a\rho,\nu^{a+1}\rho,\ldots,\nu^b\rho]\), write
\(\St(\Delta)\) for the unique irreducible quotient of
\(\nu^a\rho\times\nu^{a+1}\rho\times\cdots\times\nu^b\rho\)
\cite[\S9.1]{ZelII}. It is essentially square-integrable, and every
irreducible essentially square-integrable representation of a general
linear group arises this way \cite[Theorem~9.3]{ZelII}. We also use the
abbreviation \(\Delta(\rho,a,b)=\St(\Delta)\).

The Jacquet modules used in Section~\ref{sec:limits-smooth-necessity}
are given by the following formula \cite[Proposition~9.5]{ZelII}.

\begin{exm}\label{exm: Jac of ess disc series}
Put \(\delta=\Delta(\rho,a,b)\) and \(r=d(b-a+1)\). If \(r=s+t\)
with \(s,t>0\), then
\[
 \Jac_{P_{s,t}}^{\Gl_r}(\delta)\simeq
 \begin{cases}
 \Delta(\rho,c,b)\boxtimes\Delta(\rho,a,c-1),
   &d\mid s,\quad c=b-s/d+1,\\[3pt]
 0,&d\nmid s.
 \end{cases}
\]
\end{exm}

\subsection{Geometric Lemma for the Restriction of a Representation:}

Let \(G\) be a reductive \(p\)-adic group, \(H\subseteq G\) a closed
unimodular subgroup, and \(P=MU\) a parabolic subgroup. Let
\(\sigma\in\Rep(M)\), inflated to \(P\). Suppose that
\(P\backslash G/H\) is finite and its double cosets are locally closed.
The geometric lemma gives a filtration of
\(\restr{\Ind_P^G(\sigma)}{H}\) by these orbits
\cite[Theorem~5.2]{Bern-ZeleI}; the construction is also recalled in
\cite[\S4, equation~(4)]{OffenParabolicInduction}.

Choose representatives \(\eta_1,\ldots,\eta_N\) so that
\(Y_i=\bigcup_{j\leq i}P\eta_jH\) is open, and put
\(V_i=\{f\in\Ind_P^G(\sigma):\supp(f)\subseteq Y_i\}\).
Then \(0=V_0\subseteq V_1\subseteq\cdots\subseteq V_N=\Ind_P^G(\sigma)\)
is an \(H\)-stable filtration. The following description is the form
we shall use.

\begin{lem}[Geometric lemma]\label{lem: geom lem II}
For \(P_i=P\cap\eta_iH\eta_i^{-1}\), there is an isomorphism
\[
 \Hom_H(V_i/V_{i-1},\one_H)
 \simeq\Hom_{P_i}(\sigma,\delta_{P_i}\delta_P^{-1/2}).
\]
In particular, if \(\Ind_P^G(\sigma)\) is \(H\)-distinguished, the
space on the right is nonzero for at least one \(i\).
\end{lem}

\begin{proof}
Put \(H_i=H\cap\eta_i^{-1}P\eta_i\) and
\(\tau_i(h)=\delta_P(\eta_i h\eta_i^{-1})^{1/2}
\sigma(\eta_i h\eta_i^{-1})\). Restriction to \(P\eta_iH\), followed
by \(f\mapsto(h\mapsto f(\eta_i h))\), identifies
\(V_i/V_{i-1}\) with \(\cInd_{H_i}^H\tau_i\). Since \(H\) is
unimodular, Frobenius reciprocity gives
\(\Hom_H(\cInd_{H_i}^H\tau_i,\one_H)
\simeq\Hom_{H_i}(\tau_i,\delta_{H_i})\); see
\cite[proof of Proposition~4.1]{OffenParabolicInduction}.
Conjugating by \(\eta_i\) gives the asserted formula. Finally, a
nonzero \(H\)-invariant linear form on \(V_N\) descends to
\(V_i/V_{i-1}\) for the least \(i\) on which its restriction is nonzero.
\end{proof}

The sufficient condition for a closed orbit is stated and proved in
Proposition~\ref{prp: closed orbit contribution}; it does not require
finiteness of \(P\backslash G/H\).

The stabilizers \(P_i\) and their Levi projections for the odd model
are computed in Appendix~\ref{app:parabolic-orbits}.

\par\vspace{0.5cm}
\section{Local unfolding and the Plancherel decomposition}
\label{sec:local-unfolding}

With the notation fixed above, \(\Sp_{2n}(F)\) preserves the alternating form on
\[
W_n=\langle e_1,\ldots,e_n,f_n,\ldots,f_1\rangle_F
\]
for which \(\langle e_i,f_j\rangle=\delta_{ij}\). Put \(V=W_n\oplus Fe_0\), ordered as
\(
\{e_1,\ldots,e_n,f_n,\ldots,f_1,e_0\}
\)
and set \(G=\Gl(V)\) and \(H=\Sp(W_n)\), with \(H\) acting trivially on \(Fe_0\). Let \(N_k\) be the standard upper unitriangular subgroup of \(\Gl_k(F)\), and write
\[
\psi_k(u)=\psi\bigl(u_{1,2}+\cdots+u_{k-1,k}\bigr),\qquad u\in N_k,
\]
with \(\psi_1=1\). On \(N=N_{2n+1}\), define
\begin{equation}\label{eq:degenerate-whittaker-character}
\psi_{n,n+1}(u)=\psi\left(\sum_{i=1}^{n-1}u_{i,i+1}+\sum_{i=n+1}^{2n}u_{i,i+1}\right).
\end{equation}
Equivalently, if \(u=\left(\begin{smallmatrix}a&B\\0&d\end{smallmatrix}\right)\) relative to \(V=\langle e_1,\ldots,e_n\rangle\oplus\langle f_n,\ldots,f_1,e_0\rangle\), then \(\psi_{n,n+1}(u)=\psi_n(a)\psi_{n+1}(d)\). Likewise, on \(N_{2n}\) put \[\psi_{n,n}(u)=\psi\left(\sum_{i=1}^{n-1}u_{i,i+1}+\sum_{i=n+1}^{2n-1}u_{i,i+1}\right).\]

For a closed unimodular subgroup \(R\subset G\) and a unitary character \(\chi\) of \(R\), let \(C_c^\infty(R\backslash G,\chi)\) denote the space of locally constant functions \(f:G\to\mathbb C\) satisfying \(f(rg)=\chi(r)f(g)\) and having compact support modulo \(R\). We define \(L^2(R\backslash G,\chi)\) to be the Hilbert-space completion of \(C_c^\infty(R\backslash G,\chi)\) with respect to the norm
\[
\|f\|_2^2=\int_{R\backslash G}|f(g)|^2\,d\dot g.
\]
The group \(G\) acts by right translation. Haar measures will always be chosen compatibly with quotient integration, and additive Haar measures occurring in Fourier transforms will be self-dual with respect to \(\psi\). This is the usual realization of unitary induction:
\(L^2(R\backslash G,\chi)\simeq\uInd_R^G(\chi)\); see
\cite[Chapter~III, \S3.2, pp.~135--136]{Mackey}. We first prove the following theorem.

\begingroup
\renewcommand{\thethm}{1.1}
\begin{thm}\label{thm:local-unfolding}
There are unitary equivariant isomorphisms
\[
L^2\bigl(\Sp_{2n}(F)\backslash\Gl_{2n+1}(F)\bigr)
\simeq
L^2\bigl(N_{2n+1}\backslash\Gl_{2n+1}(F),\psi_{n,n+1}\bigr),\]and\[L^2\bigl(P(\Sp_{2n})(F)\backslash\Gl_{2n}(F)\bigr)\simeq L^2\bigl(N_{2n}\backslash\Gl_{2n}(F),\psi_{n,n}\bigr).
\]
\end{thm}
\endgroup

The proof proceeds by successive Fourier transforms, each yielding a unitary isomorphism between equivariant \(L^2\)-spaces. The following lemma and corollary provide a general framework for these intermediate isomorphisms.

Let \(W\) be a nonzero symplectic \(F\)-vector space, and set \(\cA(W)=\Sp(W)\ltimes W^\vee\), with multiplication
\[
(h,\delta)(h',\delta')=(hh',\delta h'+\delta'),
\]
where \((\delta h')(w)=\delta(h'w)\) for \(w\in W\). Let \(T\) be a second-countable, locally compact, totally disconnected group equipped with a continuous open surjective homomorphism
\[
p:T\longrightarrow\cA(W),\qquad p(t)=(h_t,\delta_t).
\]
Thus \(h_{tt'}=h_th_{t'}\) and \(\delta_{tt'}=\delta_th_{t'}+\delta_{t'}\). Fix \(\lambda_0\in W\setminus\{0\}\), and put
\[
R=\{t\in T:\delta_t=0\},\qquad
R_\star=\{t\in T:h_t\lambda_0=\lambda_0\}.
\]
Assume that \(T\), \(R\), and \(R_\star\) are unimodular. Given a smooth unitary character \(\tilde\chi\) of \(T\), define
\[
\chi=\restr{\tilde\chi}{R},\qquad
\chi_\star(r)=\tilde\chi(r)\psi\bigl(\delta_r(\lambda_0)\bigr),\quad r\in R_\star.
\]
\begin{lem}\label{lem: fourier 1}
There is a unitary \(T\)-equivariant isomorphism
\begin{equation}
U:L^2(R\backslash T,\chi)\xrightarrow{\;\sim\;}L^2(R_\star\backslash T,\chi_\star).
\end{equation}
\end{lem}

\begin{proof}
Since \(p\) is continuous, open and surjective, so is the coordinate map \(t\mapsto\delta_t\), whose fibres are the left \(R\)-cosets. It therefore induces a homeomorphism
\begin{equation}
R\backslash T\xrightarrow{\;\sim\;}W^\vee,\qquad Rt\mapsto\delta_t.
\end{equation}
Put \(q(t)=h_t^{-1}\lambda_0\). The orbit map \(\Sp(W)\to W^\times\), \(h\mapsto h^{-1}\lambda_0\), is surjective by transitivity and open because it is an analytic submersion (\(\mathfrak{sp}(W)\lambda_0=W\)). Hence \(q\) is continuous, open and surjective, and \(q(t_1)=q(t_2)\) if and only if \(t_1t_2^{-1}\in R_\star\). Consequently, \(q\) induces a homeomorphism
\begin{equation}
R_\star\backslash T\xrightarrow{\;\sim\;}W^\times,\qquad R_\star t\mapsto q(t).
\end{equation}
These homeomorphisms yield linear isomorphisms
\begin{equation}
\begin{aligned}
A &: C_c^\infty(W^\vee)\xrightarrow{\;\sim\;}C_c^\infty(R\backslash T,\chi),\\
B &: C_c^\infty(W^\times)\xrightarrow{\;\sim\;}C_c^\infty(R_\star\backslash T,\chi_\star),
\end{aligned}
\end{equation}
defined by
\[
(A\phi)(t)=\tilde\chi(t)\phi(\delta_t),\qquad
(Bv)(t)=\gamma(t)v(q(t)),
\]
where \(\gamma(t)=\tilde\chi(t)\psi\bigl(\delta_t(q(t))\bigr)\). For \(r\in R_\star\), one has
\[
q(rt)=q(t),\qquad \gamma(rt)=\chi_\star(r)\gamma(t),
\]
which gives the required left equivariance of \(Bv\). The inverses of \(A\) and \(B\) are given, for \(f\in C_c^\infty(R\backslash T,\chi)\) and \(f_\star\in C_c^\infty(R_\star\backslash T,\chi_\star)\), by
\[
\begin{aligned}
(A^{-1}f)(\delta)&=\tilde\chi(t_\delta)^{-1}f(t_\delta),&&\delta\in W^\vee,\\
(B^{-1}f_\star)(\lambda)&=\tilde\chi(s_\lambda)^{-1}\psi\bigl(-\delta_{s_\lambda}(\lambda)\bigr)f_\star(s_\lambda),&&\lambda\in W^\times,
\end{aligned}
\]
where \(t_\delta,s_\lambda\in T\) are arbitrary representatives satisfying \(\delta_{t_\delta}=\delta\) and \(q(s_\lambda)=\lambda\). These expressions are independent of the chosen representatives.

Choose mutually dual additive Haar measures \(d\delta\) on \(W^\vee\) and \(d\lambda\) on \(W\). Under the above quotient identifications, right translation by \(s\in T\) acts on \(W^\vee\) and \(W^\times\) by
\[
\delta\longmapsto\delta h_s+\delta_s,\qquad
\lambda\longmapsto h_s^{-1}\lambda,
\]
respectively. Since \(h_s\in\Sp(W)\), these actions preserve these measures. We therefore normalize the invariant quotient measures on \(R\backslash T\) and \(R_\star\backslash T\) to correspond to \(d\delta\) and \(d\lambda|_{W^\times}\). As \(|\tilde\chi|=|\gamma|=1\), the maps \(A\) and \(B\) then satisfy
\[
\|A\phi\|_{L^2(R\backslash T,\chi)}=\|\phi\|_{L^2(W^\vee)},\qquad
\|Bv\|_{L^2(R_\star\backslash T,\chi_\star)}=\|v\|_{L^2(W^\times)}.
\]

For \(V=W,W^\vee\), equip \(C_c^\infty(V)\) with the Schwartz--Bruhat topology, namely the locally convex inductive limit over fixed compact support and fixed open-lattice translation invariance. Give \(C_c^\infty(W^\times)\) the topology induced by extension by zero to \(W\). With these measures, Fourier transformation and inversion are
\begin{equation}
\widehat\phi(\lambda)=\int_{W^\vee}\phi(\delta)\psi\bigl(-\delta(\lambda)\bigr)\,d\delta,\qquad
\phi(\delta)=\int_W\widehat\phi(\lambda)\psi\bigl(\delta(\lambda)\bigr)\,d\lambda.
\end{equation}
Thus \(\cF:\phi\mapsto\widehat\phi\) is a topological isomorphism \(C_c^\infty(W^\vee)\simeq C_c^\infty(W)\) and extends to a unitary isomorphism \(L^2(W^\vee)\simeq L^2(W)\); see \cite[\S3.5 and \S\S4.4--4.5]{PoonenTate}. Set
\[
\mathcal S_0(W^\vee)
:=\left\{\phi\in C_c^\infty(W^\vee):\widehat\phi(0)=0\right\}
=\left\{\phi\in C_c^\infty(W^\vee):\int_{W^\vee}\phi(\delta)\,d\delta=0\right\},
\]
with the induced Schwartz--Bruhat topology. Since \(\widehat\phi\) is locally constant, the condition \(\widehat\phi(0)=0\) is equivalent to its vanishing on a neighbourhood of \(0\). Thus, the Fourier transformation restricts to an \(L^2\)-isometric topological isomorphism
\begin{equation}
\cF_0:\mathcal S_0(W^\vee)\xrightarrow{\;\sim\;}C_c^\infty(W^\times),\qquad
\phi\longmapsto\widehat\phi|_{W^\times},
\end{equation}
whose inverse sends \(v\) to \(\cF^{-1}\widetilde v\), where \(\widetilde v\) denotes extension by zero to \(W\). Viewed by extension by zero, \(C_c^\infty(W^\times)\) is dense in \(L^2(W)\): if \(U_j\) is a decreasing sequence of compact open neighbourhoods of \(0\) with intersection \(\{0\}\), then \(g\,1_{W\setminus U_j}\to g\) in \(L^2(W)\) for every \(g\in C_c^\infty(W)\). Fourier--Plancherel therefore gives the density of \(\mathcal S_0(W^\vee)\) in \(L^2(W^\vee)\).
Consequently, \(\cD:=A(\mathcal S_0(W^\vee))\) is dense in \(L^2(R\backslash T,\chi)\), and
\[
U=B\cF_0A^{-1}:\cD\xrightarrow{\;\sim\;}C_c^\infty(R_\star\backslash T,\chi_\star)
\]
is an isometric linear bijection between dense subspaces.

It remains to verify equivariance. Let \(\rho\) and \(\rho_\star\) denote the representations of \(T\) on \(C_c^\infty(W^\vee)\) and \(C_c^\infty(W^\times)\) obtained by transporting right translation through \(A\) and \(B\), respectively. Explicitly, for \(s\in T\),
\[
\begin{aligned}
(\rho(s)\phi)(\delta)&=\tilde\chi(s)\phi(\delta h_s+\delta_s),\\
(\rho_\star(s)v)(\lambda)&=\tilde\chi(s)\psi\bigl(\delta_s(h_s^{-1}\lambda)\bigr)v(h_s^{-1}\lambda).
\end{aligned}
\]
Indeed, \(q(ts)=h_s^{-1}q(t)\) and
\(\gamma(ts)=\gamma(t)\tilde\chi(s)\psi\bigl(\delta_s(h_s^{-1}q(t))\bigr)\).
The substitution \(\delta'=\delta h_s+\delta_s\), whose Jacobian has absolute value one, yields
\[
\widehat{\rho(s)\phi}(\lambda)
=\tilde\chi(s)\psi\bigl(\delta_s(h_s^{-1}\lambda)\bigr)\widehat\phi(h_s^{-1}\lambda).
\]
In particular, \(\mathcal S_0(W^\vee)\) is \(\rho(T)\)-stable, and
\(\cF_0\rho(s)=\rho_\star(s)\cF_0\) on this space. Hence \(U\) is a \(T\)-equivariant isometric bijection between dense subspaces and extends uniquely to the required unitary \(T\)-equivariant isomorphism.

\end{proof}
Unitary induction in stages gives the following consequence of Lemma~\ref{lem: fourier 1}.
\begin{cor}\label{cor: fourier 1}
Retain the notation and hypotheses of Lemma~\ref{lem: fourier 1}, and suppose that \(T\) is a closed subgroup of \(G\). Then there is a unitary \(G\)-equivariant isomorphism
\[
L^2(R\backslash G,\chi)\simeq L^2(R_\star\backslash G,\chi_\star).
\]
\end{cor}

\begin{proof}
The realization above, induction in stages
\cite[Chapter~III, Theorem~3.3, p.~137]{Mackey}, and
Lemma~\ref{lem: fourier 1} give
\[
\begin{aligned}
L^2(R\backslash G,\chi)
&\simeq \uInd_R^G\chi
 \simeq \uInd_T^G\bigl(\uInd_R^T\chi\bigr)\\
&\xrightarrow{\;\uInd_T^G(U)\;}
\uInd_T^G\bigl(\uInd_{R_\star}^T\chi_\star\bigr)
 \simeq \uInd_{R_\star}^G\chi_\star
 \simeq L^2(R_\star\backslash G,\chi_\star).
\end{aligned}
\]
Here \(\uInd_T^G(U)\) sends an induced section \(\Phi\) to \(g\mapsto U(\Phi(g))\); it is unitary and \(G\)-equivariant, with inverse \(\uInd_T^G(U^{-1})\), because \(U\) is unitary and \(T\)-equivariant.

\end{proof}

We now prove Theorem~\ref{thm:local-unfolding} using Lemma~\ref{lem: fourier 1} and Corollary~\ref{cor: fourier 1}.

\begin{proof}[Proof of Theorem~\ref{thm:local-unfolding}]
We first treat the odd case. We construct the desired isomorphism step by step using Corollary~\ref{cor: fourier 1}. Each application reduces the dimension of the remaining symplectic block by two. After \(n\) steps, the subgroup becomes \(N_{2n+1}\), and the resulting character is \(\psi_{n,n+1}\). It remains to identify the successive subgroups and characters.

For the first step, put
\[
U_0=\left\{u(c)=\begin{pmatrix}I_{2n}&c\\0&1\end{pmatrix}:c\in W_n\right\},\qquad T_0=HU_0.
\]
The map
\[
p_0:T_0\longrightarrow\cA(W_n),\qquad
hu(c)\longmapsto(h,\delta_c),\qquad \delta_c(w)=\langle w,c\rangle,
\]
is a topological group isomorphism. Taking \(\tilde\chi=1\) and \(\lambda_0=e_1\) in the lemma gives \(R=H\) and
\[
R_1=\Stab_H(e_1)U_0,\qquad
\chi_1(hu(c))=\psi\bigl(\langle e_1,c\rangle\bigr).
\]
The groups \(H,T_0,R_1\) are closed and unimodular. Hence, Corollary~\ref{cor: fourier 1} now gives the first unfolding isomorphism
\begin{equation}\label{eq:first-unfolding-isomorphism}
L^2(H\backslash G)\simeq L^2(R_1\backslash G,\chi_1).
\end{equation}

To continue, we describe the successive groups in block coordinates. For \(1\le j\le n\), put \(m=n-j\) and write
\[
V=E_j\oplus W_m\oplus F_{j+1},\qquad
\begin{aligned}
E_j&=\langle e_1,\ldots,e_j\rangle,\\
W_m&=\langle e_{j+1},\ldots,e_n,f_n,\ldots,f_{j+1}\rangle,\\
F_{j+1}&=\langle f_j,\ldots,f_1,e_0\rangle.
\end{aligned}
\]
Let \(T_j\) consist of the matrices
\begin{equation}\label{eq:Rj-general-matrix}
t=\begin{pmatrix}a&A&B\\0&h&C\\0&0&d\end{pmatrix},\qquad
 a\in N_j,\quad h\in\Sp(W_m),\quad d\in N_{j+1},
\end{equation}
with unrestricted off-diagonal blocks. Let \(\alpha\) be the last row of \(A\), and \(\xi\) the first column of \(C\). Define \(R_j\subset T_j\) by
\begin{equation}\label{eq:Rj-linking-relation}
\alpha={}^t\xi J_{2m}h,
\end{equation}
and put
\begin{equation}\label{eq:Rj-character}
\tilde\chi_j(t)=\psi_j(a)\psi_{j+1}(d),\qquad
\chi_j=\tilde\chi_j|_{R_j}.
\end{equation}
Here \(J_{2m}\) represents the restricted symplectic form; when \(m=0\), the middle block is absent and the relation is vacuous. For \(j=1\), these definitions recover \(R_1\) and \(\chi_1\) above.

Suppose \(j<n\). The discrepancy
\[
\delta_j(t)=\alpha-{}^t\xi J_{2m}h
\]
satisfies \(\delta_j(tt')=\delta_j(t)h_{t'}+\delta_j(t')\), by block multiplication. Hence
\[
p_j:T_j\longrightarrow\cA(W_m),\qquad t\longmapsto(h,\delta_j(t)),
\]
is a homomorphism. Replacing \(\alpha\) by \(\delta_j(t)\) is an invertible continuous change of block coordinates, so \(p_j\) is continuous, open and surjective. Moreover, \(R_j=\{t:\delta_j(t)=0\}\), and \(\tilde\chi_j\) is a smooth unitary character extending \(\chi_j\). Both \(T_j\) and \(R_j\) are closed semidirect products of unipotent groups with \(\Sp(W_m)\). They are unimodular because the determinant of the conjugation action on the unipotent radical is an algebraic character of the symplectic group and hence is trivial.
Apply the lemma with \(\lambda_0=e_{j+1}\). Its new subgroup is \(\{t\in T_j:he_{j+1}=e_{j+1}\}\). To identify it, write
\[
W_m=Fe_{j+1}\oplus W_{m-1}\oplus Ff_{j+1}.
\]
The condition \(he_{j+1}=e_{j+1}\) is equivalent to
\[
h=\begin{pmatrix}
1&{}^t\eta J_{2m-2}h'&z\\
0&h'&\eta\\
0&0&1
\end{pmatrix},\qquad h'\in\Sp(W_{m-1}).
\]
Regrouping the blocks as \(E_{j+1}\oplus W_{m-1}\oplus F_{j+2}\) therefore gives exactly \(R_{j+1}\). For the character, let \(x=\alpha(e_{j+1})\) and let \(y\) be the \(f_{j+1}\)-coordinate of \(\xi\). On this stabilizer,
\[
\delta_j(t)(e_{j+1})=x-\langle\xi,e_{j+1}\rangle=x+y.
\]
These are precisely the two new adjacent entries in the outer unipotent blocks \(a'\in N_{j+1}\) and \(d'\in N_{j+2}\). Consequently,
\[
\tilde\chi_j(t)\psi\bigl(\delta_j(t)(e_{j+1})\bigr)
=\psi_{j+1}(a')\psi_{j+2}(d')=\chi_{j+1}(t).
\]
All the hypotheses of Corollary~\ref{cor: fourier 1} are thus satisfied, giving
\begin{equation}\label{eq:successive-unfolding-isomorphism}
L^2(R_j\backslash G,\chi_j)\simeq
L^2(R_{j+1}\backslash G,\chi_{j+1}),\qquad 1\le j<n.
\end{equation}

At the last stage the symplectic block disappears:
\[
R_n=\left\{\begin{pmatrix}a&B\\0&d\end{pmatrix}:
 a\in N_n,\ d\in N_{n+1},\ B\in\Mat_{n,n+1}(F)\right\}=N_{2n+1},
\qquad \chi_n=\psi_{n,n+1}.
\]
Composing \eqref{eq:first-unfolding-isomorphism} with \eqref{eq:successive-unfolding-isomorphism} proves the first assertion. Every isomorphism is unitary and \(G\)-equivariant by Corollary~\ref{cor: fourier 1}.

For the even case, omit \(Fe_0\) and start at \(j=1\) with \(R_1=P(\Sp_{2n})\) and \(\chi_1=1\). Replace \(F_{j+1}\) by \(F_j=\langle f_j,\ldots,f_1\rangle\), take both outer blocks in \(N_j\), and set \(\tilde\chi_j(t)=\psi_j(a)\psi_j(d)\). The same maps \(p_j\) and stabilizer calculation give the successive unitary isomorphisms, now ending at \((R_n,\chi_n)=(N_{2n},\psi_{n,n})\). This proves the second assertion.

\end{proof}

\subsection{The Plancherel decompositions}

We now derive the local analytic Plancherel decompositions associated with the
Ben-Zvi--Sakellaridis--Venkatesh duality.  For \(r\geq 1\), let
\(\Pi_{\mathrm{temp}}(\Gl_r)\) denote the unitary tempered dual of \(\Gl_r(F)\), and
let \(d\mu_r\) be its Whittaker--Plancherel measure, with the
normalization determined by the Haar measures fixed above.  If
\(\pi\in\Pi_{\mathrm{temp}}(\Gl_r)\), its Hilbert space will be denoted by
\(\cH_\pi\).

\begingroup
\renewcommand{\thethm}{1.2}
\begin{thm}\label{thm:local-bzsv}
As unitary representations of their respective ambient general linear groups, one has
\begin{equation}\label{eq:bzsv-plancherel-odd}
L^2\bigl(\Sp_{2n}(F)\backslash G\bigr)
\simeq
\int_{\Pi_{\mathrm{temp}}(\Gl_n)\times\Pi_{\mathrm{temp}}(\Gl_{n+1})}^{\oplus}
\cH_{\pi_n\times\pi_{n+1}}\,
d\mu_n(\pi_n)\,d\mu_{n+1}(\pi_{n+1}).
\end{equation}
\begin{equation}\label{eq:bzsv-plancherel-even}L^2\bigl(P(\Sp_{2n})(F)\backslash \Gl_{2n}(F)\bigr)\simeq\int_{\Pi_{\mathrm{temp}}(\Gl_n)\times      \Pi_{\mathrm{temp}}(\Gl_n)}^{\oplus}\cH_{\pi_1\times\pi_2}\,d\mu_n(\pi_1)\,d\mu_n(\pi_2).\end{equation}Equivalently, if \(\phi_r\) denotes a tempered Langlands parameter for
\(\Gl_r(F)\), then the fiber indexed by \((\phi_n,\phi_{n+1})\) is the
irreducible tempered representation with parameter
\(\phi_n\oplus\phi_{n+1}\). In the even formula, the fiber indexed by\((\phi_1,\phi_2)\) is the irreducible tempered representation with parameter \(\phi_1\oplus\phi_2\).
\end{thm}
\endgroup

\begin{proof}
Let \(P=P_{n,n+1}=MU\), where \(M=\Gl_n(F)\times\Gl_{n+1}(F)\), and put
\(N_M=N_n\times N_{n+1}\) and \(\psi_M=\psi_n\boxtimes\psi_{n+1}\).
Since \(N=N_MU\) and \(\psi_{n,n+1}\) is trivial on \(U\), unitary induction in stages \cite[Chapter~III, Theorem~3.3]{Mackey} gives
\begin{equation}\label{eq:unitary-induction-in-stages}
L^2(N\backslash G,\psi_{n,n+1})
\simeq \Ind_P^G L^2(N_M\backslash M,\psi_M),
\end{equation}
where parabolic induction is normalized. We recall the Whittaker--Plancherel theorem in the form needed here. Write \(\Gl_r(F)=N_rA_rK_r\), where \(A_r\) is the diagonal subgroup and \(K_r=\Gl_r(\mathcal O_F)\), with \(\mathcal O_F\) the ring of integers of \(F\). For \(a=\operatorname{diag}(a_1,\ldots,a_r)\), put
\[
 \ell(a)=\sum_{i=1}^r\bigl|\log|a_i|_F\bigr|,
 \qquad \delta_r(a)=\prod_{i<j}|a_i/a_j|_F.
\]
The Whittaker-Schwartz space $\cC_r$, as defined in \cite[\S3.3, (3.11)]{DelormeWhittakerPlancherel}, consists of functions \(f:\Gl_r(F)\to\mathbb C\) that are right invariant under some compact open subgroup, satisfy \(f(ug)=\psi_r(u)f(g)\) for \(u\in N_r\), and obey
\[
 \nu_d(f):=\sup_{a\in A_r,\,k\in K_r}
 \delta_r(a)^{-1/2}(1+\ell(a))^d|f(ak)|<\infty
 \qquad(d\in\mathbb Z_{\geq0}).
\]
On each subspace fixed by a compact open subgroup, these seminorms define the topology; the full space carries the inductive-limit topology.

The Fourier--Whittaker transform is an equivariant isometry from
\(\cC_r\) onto the space of smooth spectral sections satisfying the
intertwining relations \cite[Proposition~5.10(i), \S7.3 and Theorem~7.7]{DelormeWhittakerPlancherel}.
Every irreducible tempered representation of \(\Gl_r(F)\) is generic, and
\(\operatorname{Wh}_{\psi_r}(\pi):=\Hom_{N_r}(\pi,\psi_r)\) is
one-dimensional \cite[Theorem~9.7]{ZelII}. Moreover, unitary induction from
square-integrable data is irreducible, and the inducing data are unique up
to conjugacy \cite[\S2, after (12)]{ABPSInnerFormsSLn}.
Thus the intertwining relations identify the spectral sections with sections
over \(\Pi_{\mathrm{temp}}(\Gl_r)\), with fiber
\(\operatorname{Wh}_{\psi_r}(\pi)\otimes\cH_\pi\). Choose measurable unit
vectors in these one-dimensional multiplicity spaces, using the Hermitian
forms in \cite[\S7.3, (7.3)]{DelormeWhittakerPlancherel}, and denote the
resulting measure by \(d\mu_r\). Restriction of the transform gives a
\(\Gl_r(F)\)-equivariant map
\[
 \mathcal F_r:C_c^\infty(N_r\backslash\Gl_r(F),\psi_r)\longrightarrow
 \int_{\Pi_{\mathrm{temp}}(\Gl_r)}^{\oplus}
       \cH_\pi\,d\mu_r(\pi)
\]
satisfying, for every \(f\) in its domain,
\[
 \int_{N_r\backslash\Gl_r(F)}|f(g)|^2\,d\dot g
 =\int_{\Pi_{\mathrm{temp}}(\Gl_r)}
       \|\mathcal F_r f(\pi)\|_{\cH_\pi}^2\,d\mu_r(\pi).
\]
Compactly supported smooth Whittaker functions are dense in \(\cC_r\)
\cite[\S3.3, (3.12)]{DelormeWhittakerPlancherel}, and the inclusion
\(\cC_r\hookrightarrow L^2(N_r\backslash\Gl_r(F),\psi_r)\) is continuous
by the estimate in \cite[\S3.3, (3.14)]{DelormeWhittakerPlancherel}.
The full transform is onto the smooth spectral sections, which are dense
in their Hilbert completion. The norm identity therefore implies that the
restricted transform also has dense image. Its domain is dense in \(L^2(N_r\backslash\Gl_r(F),\psi_r)\), so it extends by completion to a unitary \(\Gl_r(F)\)-equivariant isomorphism onto the displayed direct integral. Taking the Hilbert tensor product for \(r=n\) and \(r=n+1\), we obtain
\begin{equation}\label{eq:whittaker-plancherel-levi}
L^2(N_M\backslash M,\psi_M)\simeq
\int_{\Pi_{\mathrm{temp}}(\Gl_n)\times\Pi_{\mathrm{temp}}(\Gl_{n+1})}^{\oplus}
\cH_{\pi_n\boxtimes\pi_{n+1}}\,
d\mu_n(\pi_n)\,d\mu_{n+1}(\pi_{n+1}).
\end{equation}
Unitary induction commutes with direct integrals
\cite[Chapter~III, proposition following Theorem~3.3, pp.~137--138]{Mackey}.
Combining this fact with
Theorem~\ref{thm:local-unfolding} proves \eqref{eq:bzsv-plancherel-odd}.
The same argument for \(P_{n,n}\subset\Gl_{2n}(F)\), using the second isomorphism of
Theorem~\ref{thm:local-unfolding}, proves \eqref{eq:bzsv-plancherel-even}.

The induced representations in both formulas are irreducible and tempered by the
classification for general linear groups \cite{Bern-ZeleI,ZelII}.
Compatibility of the local Langlands correspondence with normalized parabolic induction
identifies their parameters with the indicated direct sums.  In the odd case, this
block-diagonal parameter map agrees with the prediction of the duality tables in
\cite{BenZviSakellaridisVenkatesh}.
\end{proof}

\begin{rmk}
Theorem~\ref{thm:local-unfolding} is an explicit verification, for the two models considered here, of the local unfolding procedure in \cite[Theorem~9.5.9 and \S9.5.12]{SakellaridisVenkatesh}. In the even case its conclusion is already stated in \cite[Example~9.5.17]{SakellaridisVenkatesh}, in the equivalent form of unitary parabolic induction from a product of Whittaker models. The proof above specifies the successive stabilizers and characters, verifies the unitary maps on dense subspaces and their completions, and carries out the corresponding argument for the odd model.
\end{rmk}

We now turn to smooth distinction for \(X_n^{\text{even}}\) and \(X_n^{\text{odd}}\). In the next section, we prove a sufficient condition for distinction in both models.
\par\vspace{0.5cm}
\section{Sufficient conditions for smooth distinction}
\label{sec:sufficient-distinction}

We prove Theorem~\ref{thm:intro-sufficient} using the identity orbits for \(P_{n,n+1}\) in the odd model
and \(P_{n,n}\) in the even model. We first
prove that these orbits are closed.
We then compute the relative modular characters of the closed-orbit
stabilizers and construct nonzero equivariant functionals using
Rankin--Selberg integrals. Restriction to a closed orbit then gives
distinction of the induced representations, by the argument of
\cite[Lemma~4.3]{OffenParabolicInduction}; we give the form needed here
in Proposition~\ref{prp: closed orbit contribution}.

The complete decomposition of
\(P_{l,2n+1-l}(F)\backslash\Gl_{2n+1}(F)/\Sp_{2n}(F)\), together
with the orbit representatives, stabilizers and Levi projections,
is recorded in Appendix~\ref{app:parabolic-orbits}. We use its specialization
to \(P_{1,4}\subset\Gl_5\) in the proof of Theorem~\ref{thm:intro-GL5-necessary}.

\subsection{Closed orbits}
\label{subsec:open-closed-orbits}

Put \(V=F^{2n}\) and
\(E_j=\Span\{e_1,\ldots,e_j\}\) for \(1\leq j\leq n\).

\begin{lem}\label{lem:closed-identity-orbits}
The identity orbit \(P\backslash PH\) is closed in \(P\backslash G\)
for each of the following triples:
\begin{enumerate}[label=\textup{(\arabic*)}]
\item \(G=\Gl_{2n+1}(F)\), \(H=\Sp_{2n}(F)\), and
\(P=P_{l,2n+1-l}\), where \(1\leq l\leq n\);
\item \(G=\Gl_{2n}(F)\), \(H=P(\Sp_{2n})\), and \(P=P_{n,n}\).
\end{enumerate}
\end{lem}

\begin{proof}
Since \(G\to P\backslash G\) is a quotient map, it suffices to prove
that \(PH\) is closed in \(G\). The parabolic \(P\) is the stabilizer
of \(E_d\), where \(d=l\) in \textup{(1)} and \(d=n\) in
\textup{(2)}. Thus
\[
 g\in PH\quad\Longleftrightarrow\quad
 g^{-1}E_d=hE_d\text{ for some }h\in H,
\]
since the equality on the right is equivalent to \(gh\in P\).

In \textup{(1)}, every subspace \(hE_l\), with \(h\in H\), lies in
\(V\) and is totally isotropic. Conversely, let \(W\subset V\) be
an \(l\)-dimensional totally isotropic subspace. Choose a basis
\(w_1,\ldots,w_l\) of \(W\). The map \(e_i\mapsto w_i\) is an
isometry \(E_l\to W\), since both restricted forms are zero.
Witt's extension theorem gives \(h\in\Sp_{2n}(F)\) with
\(hE_l=W\). Consequently,
\[
 g\in PH\quad\Longleftrightarrow\quad
 g^{-1}E_l\subset V\ \text{ and }\ g^{-1}E_l\text{ is totally isotropic}.
\]
To express these conditions as equations, write the first \(l\)
columns of \(g^{-1}\) as
\[
 (g^{-1}e_1,\ldots,g^{-1}e_l)
   =\begin{pmatrix}B(g)\\c(g)\end{pmatrix},
 \qquad B(g)\in\Mat_{2n,l}(F),\quad c(g)\in\Mat_{1,l}(F).
\]
Then
\[
 PH=\{g\in G:c(g)=0,\ {}^tB(g)J_{2n}B(g)=0\},
\]
which is closed because matrix inversion and multiplication are continuous.

In \textup{(2)}, every \(hE_n\), with \(h\in H\), is totally
isotropic and contains \(e_1\), since \(he_1=e_1\). Conversely,
for an \(n\)-dimensional totally isotropic subspace \(W\) containing
\(e_1\), choose a basis \(e_1,w_2,\ldots,w_n\). The isometry
\(E_n\to W\) given by \(e_1\mapsto e_1\) and \(e_i\mapsto w_i\)
for \(i\geq2\) extends by Witt's theorem to \(h\in\Sp_{2n}(F)\).
This extension fixes \(e_1\), so \(h\in H\) and \(hE_n=W\).
Putting \(A(g)=(g^{-1}e_1,\ldots,g^{-1}e_n)\), we therefore obtain
\[
 PH=\{g\in G:{}^tA(g)J_{2n}A(g)=0,\ ge_1\in E_n\}.
\]
Here \(e_1\in g^{-1}E_n\) is equivalent to \(ge_1\in E_n\),
which says that the last \(n\) entries of the first column of \(g\)
vanish. These are again closed conditions, proving \textup{(2)}.
\end{proof}

Both subgroups \(H\) in Lemma~\ref{lem:closed-identity-orbits} are
unimodular. For the mirabolic symplectic subgroup, this follows from
\(\Sp_{2n}(F)/P(\Sp_{2n})\simeq V\setminus\{0\}\), which carries the
invariant measure induced from \(V\).

\subsection{Periods on the closed orbits}
\label{subsec:closed-orbit-periods}

We now construct the equivariant functionals on the inducing data
attached to the closed orbits of Lemma~\ref{lem:closed-identity-orbits}.
The relative modular characters are computed explicitly so that the
resulting Hom-spaces agree with those in Lemma~\ref{lem: geom lem II}.

\begin{prp}\label{prp:odd-closed-orbit-period}
Put \(G=\Gl_{2n+1}(F)\), \(H=\Sp_{2n}(F)\),
\(P=P_{n,n+1}\), and \(P_H=P\cap H\).
If \(\sigma_1\) and \(\sigma_2\) are irreducible tempered representations
of \(\Gl_n(F)\) and \(\Gl_{n+1}(F)\), respectively, then
\[
 \Hom_{P_H}\left(\sigma_1\boxtimes\sigma_2,
       \delta_{P_H}\delta_P^{-1/2}\right)\neq0.
\]
\end{prp}

\begin{proof}
Put \(\sigma=\sigma_1\boxtimes\sigma_2\), inflated to \(P\).
The symplectic condition gives the Siegel stabilizer
\[
P_H
=
\left\{
\begin{pmatrix}
g&g\w_{0,n}S&0\\
0&g^\#&0\\
0&0&1
\end{pmatrix}
:
g\in\Gl_n(F),\ {}^tS=S
\right\}.
\]
Its Levi projection into
\(
M_{n,n+1}\simeq\Gl_n(F)\times\Gl_{n+1}(F)
\)
is
\[
g\longmapsto
\left(g,\diag(g^\#,1)\right).
\]
Since \(P_H\) is the Siegel parabolic subgroup of \(H\), its modular
character is
\(
\delta_{P_H}(g)=\nu(g)^{n+1}.
\)
On the other hand,
\(
\delta_P\left(g,\diag(g^\#,1)\right)
=\nu(g)^{2n+1}.
\)
Consequently,
\begin{equation}\label{eq:odd-relative-modulus}
\left.
\delta_{P_H}\delta_P^{-1/2}
\right|_{\Gl_n}
=
\nu^{1/2}.
\end{equation}

Since \(\sigma_1\) and \(\sigma_2\) are tempered, they are generic.
Choose
\(
W_1\in\cW(\sigma_1,\psi),
\qquad
W_2\in\cW(\sigma_2,\psi).
\)
For \(\Re(s)\ll0\), consider the Rankin--Selberg integral
\cite{MR701565}
\[
Z(s;W_1,W_2)
=
\int_{N_n\backslash\Gl_n(F)}
W_1(g)\,
W_2\left(\diag(g^\#,1)\right)
\nu(g)^{s-\frac12}\,dg.
\]
Since
\(
\psi(u^\#)=\psi(u)^{-1},
\ u\in N_n,
\)
the integrand is left \(N_n\)-invariant. After the change of variables \(g\mapsto g^\#\), and using
\(\sigma_1^\#\simeq\widetilde{\sigma}_1\), this is the standard local
Rankin--Selberg integral for
\(\widetilde{\sigma}_1\times\sigma_2\) at the parameter \(1-s\).
It therefore admits meromorphic continuation and is not identically zero.
The common denominator of the local Rankin--Selberg integrals gives a
uniform bound on the pole order at \(s=0\).

For \(g_0\in\Gl_n(F)\), a change of variables \(g\mapsto gg_0\) gives
\[
\begin{aligned}
&Z\left(
s;\sigma_1(g_0)W_1,
\sigma_2\left(\diag(g_0^\#,1)\right)W_2
\right)\\
&\hspace{4cm}
=
\nu(g_0)^{-s+\frac12}Z(s;W_1,W_2).
\end{aligned}
\]
Consider the Laurent expansion at \(s=0\) of the meromorphically
continued family of bilinear forms \(Z(s;\,\cdot\,,\,\cdot\,)\), and let
\(L\) be its first coefficient which is not identically zero. Comparing leading coefficients in the
preceding covariance identity gives
\[
L\left(
\sigma_1(g_0)W_1\otimes
\sigma_2\left(\diag(g_0^\#,1)\right)W_2
\right)
=
\nu(g_0)^{1/2}L(W_1\otimes W_2).
\]
The unipotent radical of \(P_H\) is contained in the unipotent radical
of \(P\) and therefore acts trivially on \(\sigma\). In view of
\eqref{eq:odd-relative-modulus}, it follows that
\[
0\neq L\in
\Hom_{P_H}
\left(
\sigma,\delta_{P_H}\delta_P^{-1/2}
\right).
\]
\end{proof}

\begin{prp}\label{prp:even-closed-orbit-period}
Put \(G=\Gl_{2n}(F)\), \(H=P(\Sp_{2n})\),
\(P=P_{n,n}\), and \(P_H=P\cap H\).
If \(\tau_1\) and \(\tau_2\) are irreducible tempered representations
of \(\Gl_n(F)\), then
\[
 \Hom_{P_H}\left(\tau_1\boxtimes\tau_2,
       \delta_{P_H}\delta_P^{-1/2}\right)\neq0.
\]
\end{prp}

\begin{proof}
Put \(\tau=\tau_1\boxtimes\tau_2\), inflated to \(P\).
The symplectic condition, together with the requirement that \(e_1\)
be fixed, gives
\[
P_H
=
\left\{
\begin{pmatrix}
g&g\w_{0,n}S\\
0&g^\#
\end{pmatrix}
:
g\in\Gl_n(F),\ ge_1=e_1,\ {}^tS=S
\right\}.
\]
Put \(q=g^\#\). Then \(q\in Q_n\), and every element of \(Q_n\) arises
in this way. Thus the Levi projection of \(P_H\) into
\(
M_{n,n}\simeq\Gl_n(F)\times\Gl_n(F)
\)
is parametrized by
\[
q\longmapsto (q^\#,q),
\qquad q\in Q_n,
\]
where, when \(n=1\), we use the convention \(Q_1=\{1\}\).

A direct computation of the modular characters gives
\begin{equation}\label{eq:even-relative-modulus}
\delta_{P_H}\bigl(\delta_P|_{P_H}\bigr)^{-1/2}
=
\one_{P_H}.
\end{equation}

The kernel of the Levi projection \(P_H\rightarrow M_{n,n}\) is
contained in the unipotent radical of \(P\), and hence acts trivially
on \(\tau\). Therefore,
\[
\Hom_{P_H}(\tau,\one_{P_H})
\simeq
\Hom_{Q_n}
\left(
\tau_1^\#\otimes\tau_2,\one_{Q_n}
\right).
\]

Since \(\tau_1\) and \(\tau_2\) are tempered, they are unitary and
generic. For \(W_1\in\cW(\tau_1,\psi)\), set
\(
W_1^\#(g)=W_1(g^\#)
\).
Then
\(
W_1^\#\in\cW(\widetilde{\tau}_1,\psi^{-1})
\).
The bilinear form
\[
L(W_1\otimes W_2)
=
\int_{N_n\backslash Q_n}
W_1(q^\#)W_2(q)\,dq,
\qquad
W_2\in\cW(\tau_2,\psi),
\]
is therefore the Rankin--Selberg period associated with
\(\widetilde{\tau}_1\otimes\tau_2\). The integral is absolutely convergent and
defines a \(Q_n\)-invariant bilinear form
\cite[Proposition~3.4]{Jo-Local-Periods}. Moreover, it is not identically
zero \cite[Theorem~3.5]{Jo-Local-Periods}. For \(n=1\), the same assertions follow directly from \(Q_1=N_1=\{1\}\).

For \(q_0\in Q_n\), right invariance of the quotient measure,
together with \((qq_0)^\#=q^\#q_0^\#\), gives
\[
L\left(
\tau_1(q_0^\#)W_1\otimes\tau_2(q_0)W_2
\right)
=
L(W_1\otimes W_2).
\]
Consequently,
\[
0\neq L\in
\Hom_{Q_n}
\left(
\tau_1^\#\otimes\tau_2,\one_{Q_n}
\right)
\simeq
\Hom_{P_H}(\tau,\one_{P_H}).
\]
In view of \eqref{eq:even-relative-modulus}, this is precisely
\(
0\neq L\in
\Hom_{P_H}
\left(
\tau,\delta_{P_H}\delta_P^{-1/2}
\right).
\)
\end{proof}

\subsection{Distinction of the induced representations}
\label{subsec:closed-orbit-criterion}

We now pass from the closed-orbit functionals constructed above to
distinction of the full induced representations.
The argument in \cite[Lemma~4.3]{OffenParabolicInduction} is stated for
symmetric spaces. The following proposition records the same restriction
and Frobenius-reciprocity argument for the closed unimodular subgroups
considered here.

\begin{prp}\label{prp: closed orbit contribution}
Let \(G\) be a reductive \(p\)-adic group, let \(H\) be a closed unimodular
subgroup of \(G\), and let \(P=MU\) be a parabolic subgroup.
Assume that the identity orbit \(P\backslash PH\) is closed in
\(P\backslash G\), and put \(P_0=P\cap H\).
If \(\sigma\in\Rep(M)\), inflated to \(P\), satisfies
\[
 \Hom_{P_0}\left(\sigma,
       \delta_{P_0}\delta_P^{-1/2}\right)\neq0,
\]
then \(\Ind_P^G(\sigma)\) is \(H\)-distinguished.
\end{prp}

\begin{proof}
Restriction to the closed orbit is a surjective \(H\)-equivariant map
\[
 r:\Ind_P^G(\sigma)\longrightarrow
 \mathcal C_\sigma
 :=\cInd_{P_0}^H\bigl((\delta_P^{1/2}\sigma)|_{P_0}\bigr).
\]
Here \(P_0\backslash H\) is compact, and surjectivity follows from
extension of locally constant sections from a closed subset of
\(P\backslash G\). The Frobenius-reciprocity identification underlying
Geometric Lemma~\ref{lem: geom lem II} gives
\[
 \Hom_H(\mathcal C_\sigma,\one_H)
 \simeq
 \Hom_{P_0}\left(\sigma,
       \delta_{P_0}\delta_P^{-1/2}\right).
\]
If \(\ell\) is a nonzero element of the left-hand side, then
\(\ell\circ r\) is a nonzero \(H\)-invariant functional on
\(\Ind_P^G(\sigma)\).
\end{proof}

\begin{restatedthm}{thm:intro-sufficient}
Let \(n\geq1\).
\begin{enumerate}[label=\textup{(\arabic*)},leftmargin=2em]
\item If \(\sigma_1\) and \(\sigma_2\) are irreducible tempered representations of \(\GL_n(F)\) and \(\GL_{n+1}(F)\), respectively, then \(\sigma_1\times\sigma_2\) is \(\Sp_{2n}(F)\)-distinguished.
\item If \(\tau_1\) and \(\tau_2\) are irreducible tempered representations of \(\GL_n(F)\), then \(\tau_1\times\tau_2\) is \(P(\Sp_{2n})\)-distinguished.
\end{enumerate}
\end{restatedthm}

\begin{proof}
For \textup{(1)}, take
\[
 (G,H,P)=
 (\Gl_{2n+1}(F),\Sp_{2n}(F),P_{n,n+1}),
 \qquad \sigma=\sigma_1\boxtimes\sigma_2.
\]
The identity orbit is closed by
Lemma~\ref{lem:closed-identity-orbits}\textup{(1)}, and
Proposition~\ref{prp:odd-closed-orbit-period} supplies the required
nonzero equivariant functional. Proposition~\ref{prp: closed orbit contribution}
therefore gives \(\Sp_{2n}(F)\)-distinction of
\(\sigma_1\times\sigma_2\).

For \textup{(2)}, take
\[
 (G,H,P)=(\Gl_{2n}(F),P(\Sp_{2n}),P_{n,n}),
 \qquad \sigma=\tau_1\boxtimes\tau_2.
\]
The subgroup \(H\) is unimodular, as noted in
Subsection~\ref{subsec:open-closed-orbits}.
Lemma~\ref{lem:closed-identity-orbits}\textup{(2)} and
Proposition~\ref{prp:even-closed-orbit-period} verify the remaining
hypotheses of Proposition~\ref{prp: closed orbit contribution}, which
proves \(P(\Sp_{2n})\)-distinction of \(\tau_1\times\tau_2\).
\end{proof}

\par\vspace{0.5cm}
\section{Smooth distinction for \texorpdfstring{$X_n^{\mathrm{odd}}$ and $X_n^{\mathrm{even}}$}{the odd and even models}}
\label{sec:limits-smooth-necessity}

We now ask whether the sufficient conditions of the previous section are also necessary for smooth distinction of irreducible tempered representations. Although we do not settle the question in full generality, we obtain precise results for several classes of representations. More precisely, we prove that essentially square-integrable representations
of \(\Gl_{2n+1}(F)\) are not \(\Sp_{2n}(F)\)-distinguished, and that
supercuspidal representations of \(\Gl_{2n}(F)\) are not
\(P(\Sp_{2n})\)-distinguished. For the latter model, we construct
square-integrable counterexamples on \(\Gl_2(F)\) and \(\Gl_4(F)\)
lying outside the corresponding Plancherel support. We then prove that
an irreducible tempered representation of \(\Gl_5(F)\) is
\(\Sp_4(F)\)-distinguished precisely when it is induced from tempered
representations on \(\Gl_2(F)\times\Gl_3(F)\). In particular, the
\(\Gl_4(F)\) counterexample does not produce an additional distinguished
representation of \(\Gl_5(F)\). Thus, it remains open whether the sufficient condition is also necessary
for \(X_n^{\mathrm{odd}}\) in general.

\subsection{Essentially square-integrable representations of \texorpdfstring{\(\Gl_{2n+1}(F)\)}{GL(2n+1,F)}}

Put
\[
 \widetilde G=\Gl_{2n+1}(F),\qquad G=\Gl_{2n}(F),\qquad
 H=\Sp_{2n}(F)\subset G,
\]
where $G$ is embedded in $\widetilde G$ by
$g\mapsto\operatorname{diag}(g,1)$.
Let $N=N_{2n}$ be the standard maximal unipotent subgroup of $G$, and fix the
nondegenerate character
\[
 \psi_N(u)=\psi(u_{1,2}+\cdots+u_{2n-1,2n}).
\]

\begin{restatedthm}{thm:intro-discrete-vanishing}
If \(\pi\) is an irreducible, essentially square-integrable representation of \(\GL_{2n+1}(F)\), then
\[
\Hom_{\Sp_{2n}(F)}(\pi,\mathbf1)=0.
\]
\end{restatedthm}

\begin{proof}
We first treat an irreducible supercuspidal representation $\rho$ of
$\widetilde G$. The Kirillov model gives an isomorphism of $G$-modules
\begin{equation}\label{eq:odd-kirillov-restriction}
 \restr{\rho}{G}\simeq \cInd_N^G\psi_N;
\end{equation}
see Proposition~\ref{prp:supercuspidal-restriction-mirabolic}. We claim that
\begin{equation}\label{eq:gelfand-graev-no-symplectic-period}
 \Hom_H\!\left(\cInd_N^G\psi_N,\one\right)=0.
\end{equation}
To prove this, we use the following special case of
\cite[Proposition~1]{OffSayDisj}, with \((r,2k)=(2n,0)\) and
\((r',2k')=(0,2n)\): any distribution $D$ on $G$ satisfying
\[
 D(L_uR_hf)=\psi_N(u)^{-1}D(f)
 \qquad (u\in N,\ h\in H,\ f\in C_c^\infty(G))
\]
vanishes. Here
\[
 (L_uf)(g)=f(u^{-1}g),\qquad (R_hf)(g)=f(gh).
\]
Suppose now that $\Lambda$ is a nonzero element of the Hom-space in
\eqref{eq:gelfand-graev-no-symplectic-period}. The averaging map
\[
 \mathcal A_{\psi_N}:C_c^\infty(G)\longrightarrow\cInd_N^G\psi_N,
 \qquad
 \mathcal A_{\psi_N}(f)(g)
   =\int_N\psi_N(u)^{-1}f(ug)\,du,
\]
is surjective. Indeed, for \(W\in\cInd_N^G\psi_N\), ordinary averaging
\cite[\S1.21, Lemma~(a)]{Bern-Zele-Survey} supplies
\(\beta\in C_c^\infty(G)\) such that \(\int_N\beta(ug)\,du=1\) on
the support of \(W\) modulo \(N\). Then
\(\mathcal A_{\psi_N}(\beta W)=W\). Consequently,
$D=\Lambda\circ\mathcal A_{\psi_N}$ is a nonzero distribution on $G$.
A change of variables gives
\[
 \mathcal A_{\psi_N}(L_uR_hf)
   =\psi_N(u)^{-1}R_h\mathcal A_{\psi_N}(f).
\]
Since $\Lambda$ is $H$-invariant, $D$ satisfies the equivariance condition
above and must therefore vanish, a contradiction. This proves
\eqref{eq:gelfand-graev-no-symplectic-period}, and hence
$\Hom_H(\rho,\one)=0$ by \eqref{eq:odd-kirillov-restriction}.

Now let $\delta$ be any irreducible essentially square-integrable
representation of $\widetilde G$, and choose an irreducible supercuspidal
representation $\rho$ of the same group. From \cite[Theorem~5.7]{Chan-Savin} we know that any two
irreducible essentially square-integrable representations of $\Gl_{m+1}(F)$
become isomorphic after restriction to $\Gl_m(F)$. Consequently
\[
 \restr{\delta}{G}\simeq\restr{\rho}{G},
 \qquad
 \Hom_H(\delta,\one)\simeq\Hom_H(\rho,\one)=0.
\]

\end{proof}

Theorem~\ref{thm:intro-discrete-vanishing} raises the same question for
\((\Gl_{2n}(F),P(\Sp_{2n}))\). We first exclude supercuspidal
representations, then give square-integrable examples showing that the
two models behave differently.
\subsection{Supercuspidal representations of \texorpdfstring{\(\Gl_{2n}(F)\)}{GL(2n,F)}}

We use \cite[Corollary~1 and Remark~2]{DijolsPrasad}, stated here in our
notation. Although the corollary is formulated for a unitary group,
Remark~2 explains that its argument applies to \(\Gl_{2n}\) and uses only
invariance under the Klingen mirabolic subgroup of \(\Sp_{2n}\). For
\(\boldsymbol\epsilon=(\epsilon_1,\ldots,\epsilon_{n-1})
\in\{0,1\}^{n-1}\), define a character of \(N_{2n}\) by
\[
 \psi_{\boldsymbol\epsilon}(u)=
 \psi\left(\sum_{j=1}^{n-1}\epsilon_j
 (u_{j,j+1}+u_{2n-j,\,2n-j+1})\right),
\]
using the ordered basis \(e_1,\ldots,e_n,f_n,\ldots,f_1\).
For \(n=1\), this is the trivial character. It follows that
any smooth representation \(\sigma\) of \(\Gl_{2n}(F)\) satisfies
\begin{equation}\label{eq:even-degenerate-whittaker-form}
 \Hom_{P(\Sp_{2n})}(\sigma,\one)\neq0
 \quad\Longrightarrow\quad
 \Hom_{N_{2n}}(\sigma,\psi_{\boldsymbol\epsilon})\neq0
 \quad\text{for some }\boldsymbol\epsilon.
\end{equation}

\begin{prp}\label{prp:even-supercuspidal-vanishing}
Every irreducible supercuspidal representation \(\pi\) of
\(\Gl_{2n}(F)\) satisfies
\[
 \Hom_{P(\Sp_{2n})}(\pi,\one)=0.
\]
\end{prp}

\begin{proof}
If \(\pi\) were \(P(\Sp_{2n})\)-distinguished,
\eqref{eq:even-degenerate-whittaker-form} would give a nonzero
\(\ell\in\Hom_{N_{2n}}(\pi,\psi_{\boldsymbol\epsilon})\).
The character \(\psi_{\boldsymbol\epsilon}\) does not involve the
entry \(u_{n,n+1}\), so it is trivial on
\[
 U_{n,n}=\left\{
 \begin{pmatrix}I_n&X\\0&I_n\end{pmatrix}:X\in\Mat_{n,n}(F)
 \right\}.
\]
Thus \(\ell\) factors nontrivially through \(\pi_{U_{n,n}}\), and hence
\(\Jac_{P_{n,n}}^{\Gl_{2n}}(\pi)\neq0\). This contradicts
supercuspidality.
\end{proof}

\subsection{Square-integrable counterexamples in ranks two and four}
\label{subsec:even-discrete-counterexamples}

In this subsection, we prove that Proposition~\ref{prp:even-supercuspidal-vanishing} does not extend from
supercuspidal to square-integrable representations. We begin with \(\Gl_2(F)\).

Since $\Sp_2(F)=\Sl_2(F)$, we have
\(
 P(\Sp_2)=\left\{\begin{pmatrix}1&x\\0&1\end{pmatrix}:x\in F\right\}=N_2.
\)
Let $\mu$ be a unitary character of $F^\times$ and put
$\delta_\mu=(\mu\circ\det)\otimes\St_2$.
The segment Jacquet-module formula (Example~\ref{exm: Jac of ess disc series}) gives
\[
 \Jac_{P_{1,1}}^{\Gl_2}(\delta_\mu)
 \simeq\nu^{1/2}\mu\boxtimes\nu^{-1/2}\mu.
\]
Consequently, the unnormalized coinvariant space $(\delta_\mu)_{N_2}$
is one-dimensional, and hence
\[
 \dim\Hom_{P(\Sp_2)}(\delta_\mu,\one)
 =\dim\bigl((\delta_\mu)_{N_2}\bigr)^\vee=1.
\]
On the other hand, Theorem~\ref{thm:intro-plancherel} for $n=1$
identifies the Plancherel support of $N_2\backslash\Gl_2(F)$ with the
unitary principal series $\chi_1\times\chi_2$, where $\chi_1,\chi_2$
are unitary characters of $F^\times$. The square-integrable representation
$\delta_\mu$ lies outside this support. Thus the smooth tempered
distinguished spectrum is already strictly larger than the $L^2$-spectrum
in rank two.

For \(\Gl_4(F)\), the following family gives further examples.
\begin{restatedthm}{thm:intro-even-counterexample}
Let \(\rho\) be an irreducible unitary supercuspidal representation of \(\GL_2(F)\), and put
\[
\delta_\rho
   =\operatorname{St}
      \bigl([\nu^{-1/2}\rho,\nu^{1/2}\rho]\bigr),
\]
where \(\operatorname{St}(\Delta)\) denotes the essentially square-integrable representation attached to the Zelevinsky segment \(\Delta\) \cite{ZelII}. Then \(\delta_\rho\) is an irreducible square-integrable representation of \(\GL_4(F)\), and
\[
\dim\Hom_{P(\Sp_4)}(\delta_\rho,\mathbf1)=1.
\]
\end{restatedthm}

\begin{proof}
Use the ordered symplectic basis $e_1,e_2,f_2,f_1$, and put
\[
 G=\Gl_4(F),\qquad
 H=P(\Sp_4)=\{h\in\Sp_4(F):he_1=e_1\}.
\]
Let $U=U_{1,2,1}$ be the unipotent radical of $P_{1,2,1}$. Write
\[
 u(x,y,z)=
 \begin{pmatrix}
 1&x_1&x_2&z\\
 0&1&0&y_1\\
 0&0&1&y_2\\
 0&0&0&1
 \end{pmatrix},
 \qquad
 x=(x_1,x_2),\quad y={}^t(y_1,y_2),
\]
and set
\[
 J_2=
 \begin{pmatrix}
 0&1\\
 -1&0
 \end{pmatrix}.
\]
A direct calculation gives
\[
 U':=U\cap H
   =\{u(x,J_2{}^tx,z):x\in F^2,\ z\in F\}.
\]
Moreover,
\(
 S=\{\operatorname{diag}(1,s,1):s\in\Sl_2(F)\}
   \simeq\Sl_2(F),
 \
 H=S\ltimes U'.
\).
The map
\[
 q:U\longrightarrow F^2,\qquad
 q(u(x,y,z))=y-J_2{}^tx
   ={}^t(y_1-x_2,\,y_2+x_1),
\]
is a surjective homomorphism with kernel $U'$, and the conjugation action of
$S$ on $U/U'\simeq F^2$ is the standard action.

Define
\[
 \xi(u(x,y,z)):=\psi(x_1+y_2).
\]
This character is trivial on $U'$ and, through $q$, corresponds to
\(
 \bar\xi(a_1,a_2)=\psi(a_2).
\)
The action of $S$ on the character group of $U/U'$ has two orbits, namely
the trivial character and the nontrivial characters. The stabilizer of
$\bar\xi$ is
\(
 N=\{I_4+tE_{2,3}:t\in F\}\subset S.
\)
The exact sequence supplied by smooth Clifford theory
\cite[Proposition~1]{DijolsPrasad} therefore gives, for every smooth
$G$-representation $\sigma$, an exact sequence of $S$-modules
\begin{equation}\label{eq:rank-four-clifford-sequence}
 0\longrightarrow\cInd_N^S(\sigma_{U,\xi})
 \longrightarrow\sigma_{U'}
 \longrightarrow\sigma_U
 \longrightarrow0.
\end{equation}

Let $\rho$ be an irreducible unitary supercuspidal representation of
$\Gl_2(F)$ and put
\[
 \delta_\rho
   =\St\bigl([\nu^{-1/2}\rho,\nu^{1/2}\rho]\bigr).
\]
This is an irreducible square-integrable representation of $G$, and we have \[\Hom_H(\delta_\rho, \one)\cong \Hom_S\left((\delta_\rho)_{U'}, \one\right).\]
Its Jacquet
module for $P_{3,1}$ vanishes by the Jacquet-module formula for a segment
representation \cite[Proposition~9.5]{ZelII}: the cuspidal degree of $\rho$ is
two, whereas the Levi contains a block of degree three. By transitivity of
Jacquet modules,
\(
 \Jac_{P_{1,2,1}}^G(\delta_\rho)=0.
\)

Taking $\sigma=\delta_\rho$ in
\eqref{eq:rank-four-clifford-sequence} therefore gives
\begin{equation}\label{eq:rank-four-Uprime}
 (\delta_\rho)_U=0,\qquad (\delta_\rho)_{U'}
 \simeq
 \cInd_N^S\bigl((\delta_\rho)_{U,\xi}\bigr).
\end{equation}

Let $N_4$ be the standard maximal unipotent subgroup of $G$ and define
\[
 \psi_{2,2}(v)=\psi(v_{1,2}+v_{3,4}),
 \ v\in N_4.
\]
By \eqref{eq:rank-four-Uprime} and Frobenius reciprocity, and since
\(N_4=U\rtimes N\), with
\(\psi_{2,2}|_U=\xi\) and \(\psi_{2,2}|_N=\one\), we have
\begin{equation}\label{eq:rank-four-semi-Whittaker}
 \Hom_S\left((\delta_\rho)_{U'},\one\right)
\simeq\Hom_N\bigl((\delta_\rho)_{U,\xi},\one\bigr)
 \simeq
 \Hom_{N_4}(\delta_\rho,\psi_{2,2})\simeq\mathbb C.
\end{equation}
To compute the last space, put $R=U_{2,2}$. Then
\(
 N_4=(N_2\times N_2)\ltimes R,
\)
and $\psi_{2,2}$ is trivial on $R$. Moreover, the Jacquet-module formula gives
\[
 \Jac_{P_{2,2}}^G(\delta_\rho)
 \simeq
 \nu^{1/2}\rho\boxtimes\nu^{-1/2}\rho
\]
for our standard upper triangular parabolic. Both factors are generic, and \(\delta_{P_{2,2}}^{1/2}\) is trivial on \(N_2\times N_2\). Taking Whittaker functionals on the two factors therefore identifies the last space in \eqref{eq:rank-four-semi-Whittaker} with a tensor product of two one-dimensional spaces. Hence,
\[
\dim\Hom_{P(\Sp_4)}(\delta_\rho,\one)=1.
\]

\end{proof}
\begin{rmk}\label{rmk:GL4-Steinberg-multiplicity}
Let $\rho$ be a unitary character of $F^\times$, and put
\[
 \delta=\St([\nu^{-3/2}\rho,\nu^{-1/2}\rho,\nu^{1/2}\rho,\nu^{3/2}\rho]),
 \qquad D=\St([\nu^{-1/2}\rho,\nu^{1/2}\rho]).
\]
Then
\[
\dim\Hom_{P(\Sp_4)}(\delta,\one)\leq1.
\]
Indeed, retain the notation
$U,U',S,N,\xi$ of the preceding proof. The Jacquet-module formula
\cite[Proposition~9.5]{ZelII} gives
\[
 \begin{aligned}
 \Jac_{P_{1,2,1}}^{\Gl_4}(\delta)
 &\simeq\nu^{3/2}\rho\boxtimes D\boxtimes\nu^{-3/2}\rho,\\
 \Jac_{P_{2,2}}^{\Gl_4}(\delta)
 &\simeq\nu D\boxtimes\nu^{-1}D.
 \end{aligned}
\]
The normalizing character for $P_{1,2,1}$ is trivial on $S$, and genericity of $D$ gives the vanishing below by \cite[Theorem~3.2.2]{HR90}. Applying $\Hom_S(-,\one)$ to
\eqref{eq:rank-four-clifford-sequence} with $\sigma=\delta$, and using
Frobenius reciprocity as in \eqref{eq:rank-four-semi-Whittaker}, gives
\begin{equation}\label{eq:rank-four-semi-Whittaker 2}
\begin{gathered}
 \Hom_S(\delta_U,\one)=\Hom_{\Sl_2(F)}(D,\one)=0,\\
 \Hom_{P(\Sp_4)}(\delta,\one)
 \hookrightarrow\Hom_N(\delta_{U,\xi},\one)
 \simeq\Hom_{N_4}(\delta,\psi_{2,2})\simeq\mathbb C.
\end{gathered}
\end{equation}
The last space is one-dimensional: $\psi_{2,2}$ is trivial on $U_{2,2}$,
so it is the Whittaker-functional space of $\nu D\boxtimes\nu^{-1}D$
for the nondegenerate character on $N_2\times N_2$. Both factors are
generic and have unique Whittaker functionals
\cite[Corollary~8.3 and Theorem~9.7]{ZelII}. This proves the bound.

We will use this bound in the proof of Theorem~\ref{thm:intro-GL5-necessary}; nonvanishing is not needed there.
\end{rmk}

In the next subsection, we study smooth distinction for
\((\Gl_5(F),\Sp_4(F))\).

\subsection{Tempered distinction for \texorpdfstring{\(\Gl_5(F)\)}{GL(5,F)}}
\label{subsec:GL5-necessary}

We prove Theorem~\ref{thm:intro-GL5-necessary}, which states that an irreducible tempered representation of \(\Gl_5(F)\)
can be \(\Sp_4(F)\)-distinguished only if it is induced from tempered
data on \(\Gl_2(F)\times\Gl_3(F)\), giving the converse of
Theorem~\ref{thm:intro-sufficient}. The tempered classification and
Theorem~\ref{thm:intro-discrete-vanishing} reduce the necessary condition
to showing that
\[
 \Hom_{\Sp_4(F)}(\chi\times\delta,\one)=0
\]
for every unitary character \(\chi\) of \(F^\times\) and every
irreducible square-integrable representation \(\delta\) of \(\Gl_4(F)\).

To prove this vanishing statement, we use the geometric lemma (Lemma~\ref{lem: geom lem II}) and the orbit calculations in
Appendix~\ref{app:parabolic-orbits} to reduce this assertion to the
vanishing of the invariant functionals on
\(C_c^\infty(F^4,V_\delta)\), where \(V_\delta\) is the space of
\(\delta\). Restriction to functions supported away from \(0\)
determines how any such functional must transform under scalar
dilations. Fourier transform with respect to the symplectic pairing
gives an incompatible transformation rule, forcing the functional
to vanish. This argument also applies to the representation
\(\delta_\rho\) of Theorem~\ref{thm:intro-even-counterexample}:
although the open orbit admits a nonzero invariant functional, that
functional does not extend to \(\chi\times\delta_\rho\).

\begin{restatedthm}{thm:intro-GL5-necessary}
Let \(\pi\) be an irreducible tempered representation of \(\Gl_5(F)\). Then
\[
 \Hom_{\Sp_4(F)}(\pi,\one)\neq0
 \quad\Longleftrightarrow\quad
 \pi\simeq\sigma_2\times\sigma_3
\]
for irreducible tempered representations \(\sigma_i\) of \(\Gl_i(F)\), \(i=2,3\).
\end{restatedthm}

\begin{proof}
The sufficient condition is Theorem~\ref{thm:intro-sufficient}. Every irreducible tempered representation of a general linear group is an irreducible normalized induction of unitary square-integrable representations, determined up to reordering of the factors; see \cite[\S2, after (12)]{ABPSInnerFormsSLn}. Among the partitions of \(5\), only \(5\) and \(4+1\) have no collection of parts summing to \(2\). All other partitions therefore give a presentation \(\sigma_2\times\sigma_3\) as in the statement. Theorem~\ref{thm:intro-discrete-vanishing} excludes the square-integrable case. It remains to prove
\begin{equation}\label{eq:GL5-chi-delta-vanishing}
 \Hom_H(\Pi,\one)=0,\qquad
 G=\Gl_5(F),\quad H=\Sp_4(F),\quad
 P=P_{1,4},\quad \Pi=\chi\times\delta,
\end{equation}
where \(\chi\) is unitary and \(\delta\) is an irreducible square-integrable representation of \(\Gl_4(F)\).

We first describe the orbit filtration of \(\Pi|_H\).
Put \(W=F^4\), with ordered basis \(e_1,e_2,f_2,f_1\), and embed \(H\) in \(G\) by \(h\mapsto\operatorname{diag}(h,1)\). The specialization \(n=2\), \(l=1\) of Appendix~\ref{app:parabolic-orbits}, in particular \eqref{eq:odd-orbit-decomposition}, gives the representatives
\[
 g_0^1=1,\qquad g_{1,0}=\omega,\qquad
 g_{2,0}=\omega u(e_1),\qquad
 \omega=(1,5),\quad u(x)=\begin{pmatrix}I_4&x\\0&1\end{pmatrix}.
\]
Under \(Pg\mapsto g^{-1}Fe_1\), the orbits \(P\backslash\cO_0^1\), \(P\backslash\cO_{1,0}\), and \(P\backslash\cO_{2,0}\) are, respectively,
\[
 \mathbb P(W),\qquad \{Fe_5\},\qquad
 \{F(e_5+x):x\in W\setminus\{0\}\}.
\]
Thus \(\cO_{2,0}\) is open, and \(P\backslash(\cO_{2,0}\cup\cO_{1,0})\) is the affine chart complementary to \(\mathbb P(W)\).

Use the ordered basis \((e_1,e_5,e_2,f_2,f_1)\) of \(F^5\), so that the four-dimensional Levi factor of \(P\) is identified with \(W\) through \(\omega\). This only conjugates \(\delta\) by a fixed change of basis. Put
\[
 J_2=\begin{pmatrix}0&1\\-1&0\end{pmatrix},\qquad
 J_4=\begin{pmatrix}0&0&1\\0&J_2&0\\-1&0&0\end{pmatrix}.
\]
In these coordinates, the stabilizer formulas \eqref{eq:base-stabilizer} and \eqref{eq:odd-orbit-stabilizers} give
\begin{equation}\label{eq:GL5-stabilizers}
\begin{gathered}
 P_0^1=\left\{\begin{pmatrix}
 t&0&x&z\\
 0&1&0&0\\
 0&0&s&t^{-1}sJ_2{}^tx\\
 0&0&0&t^{-1}
 \end{pmatrix}\right\},\\
 P_{1,0}=\{\operatorname{diag}(1,h):h\in\Sp_4(F)\},\qquad
 P_{2,0}=\{\operatorname{diag}(1,h):h\in P(\Sp_4)\}.
\end{gathered}
\end{equation}
Here \(t\in F^\times\), \(s\in\Sl_2(F)\), \(x\in F^{1\times2}\), and \(z\in F\); the block sizes in the first matrix are \((1,1,2,1)\). Writing \(p_0(t,s,x,z)\) for that matrix, their modulus characters satisfy
\[
 \delta_{P_0^1}(p_0(t,s,x,z))=|t|^4,\qquad
 \delta_P^{1/2}(p_0(t,s,x,z))=|t|^{5/2},\qquad
 \delta_{P_{i,0}}=\delta_P|_{P_{i,0}}=1\quad(i=1,2).
\]
As in \eqref{eq:first-cell-levi-projection}, \eqref{eq:M1a}, and \eqref{eq:M2a}, their Levi images are
\begin{equation}\label{eq:GL5-levi-projections}
\begin{aligned}
 M_0^1&=\left\{\left(t,A(t,s,y)\right):
 A(t,s,y)=\begin{pmatrix}1&0&0\\0&s&y\\0&0&t^{-1}\end{pmatrix},\quad
 t\in F^\times,\ s\in\Sl_2(F),\ y\in F^2\right\},\\
 M_{1,0}&=\{(1,h):h\in\Sp_4(F)\},\qquad
 M_{2,0}=\{(1,h):h\in P(\Sp_4)\}.
\end{aligned}
\end{equation}

Apply Lemma~\ref{lem: geom lem II} with the orbits ordered as \(\cO_{2,0},\cO_{1,0},\cO_0^1\). The filtration and quotient descriptions in its proof give
\begin{equation}\label{eq:GL5-orbit-filtration}
\begin{gathered}
 \mathcal A=\{f\in\Pi:\operatorname{supp}(f)\subset\cO_{2,0}\},\quad
 \mathcal B=\{f\in\Pi:\operatorname{supp}(f)\subset\cO_{2,0}\cup\cO_{1,0}\},\quad
 \mathcal C=\Pi|_H/\mathcal B,\\
 0\subset\mathcal A\subset\mathcal B\subset\Pi|_H,\qquad
 \mathcal A\simeq\cInd_{P(\Sp_4)}^H(\delta|_{P(\Sp_4)}),\qquad
 \mathcal B/\mathcal A\simeq\delta|_H.
\end{gathered}
\end{equation}
Here compact induction is unnormalized. Lemma~\ref{lem: geom lem II}, with the modulus characters above, gives the corresponding Hom-spaces:
\begin{equation}\label{eq:GL5-orbit-Hom}
\begin{gathered}
 \Hom_H(\mathcal A,\one)\simeq\Hom_{P(\Sp_4)}(\delta,\one),\qquad
 \Hom_H(\mathcal B/\mathcal A,\one)\simeq\Hom_{\Sp_4(F)}(\delta,\one)=0,\\
 \Hom_H(\mathcal C,\one)\simeq
 \Hom_{M_0^1}(\chi\boxtimes\delta,\eta),\qquad
 \eta(t,A(t,s,y))=|t|^{3/2}.
\end{gathered}
\end{equation}
The vanishing of \(\Hom_H(\mathcal B/\mathcal A,\one)\) follows from the genericity of \(\delta\) and \cite[Theorem~3.2.2]{HR90}.

We next prove \(\Hom_H(\mathcal C,\one)=0\). Restricting the last Hom-space in \eqref{eq:GL5-orbit-Hom} to \(t=1\) gives an injection into \(\Hom_K(\delta,\one)\), where
\[
 K=\left\{\begin{pmatrix}1&0&0\\0&s&y\\0&0&1\end{pmatrix}:
 s\in\Sl_2(F),\ y\in F^2\right\}.
\]
Indeed, both \(\chi(t)\) and \(\eta\) are trivial on this subgroup. Let
\(\Gamma_j=\cInd_{N_j}^{\Gl_j(F)}\psi_j\). As in the proof of Theorem~\ref{thm:intro-discrete-vanishing}, \cite[Theorem~5.7]{Chan-Savin} and Proposition~\ref{prp:supercuspidal-restriction-mirabolic} give
\[
 \delta|_{\operatorname{diag}(1,\Gl_3(F))}\simeq\Gamma_3.
\]
In this \(\Gl_3\), the group \(K\) is \(\Sl_2\ltimes U_{2,1}\). The Jacquet-module formula \cite[Lemma~4.5]{ChanProjectivity} gives
\[
 \Jac_{P_{2,1}}^{\Gl_3}(\Gamma_3)\simeq\Gamma_2\boxtimes\Gamma_1.
\]
The normalizing character is trivial on \(\Sl_2\). Therefore
\[
 \Hom_K(\delta,\one)
 \simeq\Hom_{\Sl_2(F)}(\Gamma_2\otimes\Gamma_1,\one)=0.
\]
For the last equality, \(\Sl_2\) acts only on \(\Gamma_2\), and
\(\Hom_{\Sl_2(F)}(\Gamma_2,\one)=0\) is exactly
\eqref{eq:gelfand-graev-no-symplectic-period} with \(n=1\), proved in Theorem~\ref{thm:intro-discrete-vanishing}. Thus \(\Hom_H(\mathcal C,\one)=0\).

Since \(\Hom_H(\mathcal C,\one)=\Hom_H(\mathcal B/\mathcal A,\one)=0\), restriction gives injections
\begin{equation}\label{eq:GL5-restriction-injections}
 \Hom_H(\Pi,\one)\hookrightarrow\Hom_H(\mathcal B,\one)
 \hookrightarrow\Hom_H(\mathcal A,\one).
\end{equation}
In particular, proving \(\Hom_H(\mathcal B,\one)=0\) proves \eqref{eq:GL5-chi-delta-vanishing}.

The assignment \(f\longmapsto\Phi_f,\ \Phi_f(x)=f(\omega u(x))\) identifies
\begin{equation}\label{eq:GL5-affine-model}
 \mathcal B\simeq C_c^\infty(W,V_\delta),\qquad
 \mathcal A\simeq C_c^\infty(W\setminus\{0\},V_\delta).
\end{equation}
Conversely, \(\Phi\) determines \(f\) on this chart by the defining \(P\)-equivariance of \(\Pi\), and extension by zero gives a smooth section because \(\Phi\) has compact support inside the open chart. In the Levi coordinates fixed above,
\[
 u(x)\operatorname{diag}(h,1)
 =\operatorname{diag}(h,1)u(h^{-1}x),\qquad
 \omega\operatorname{diag}(h,1)\omega^{-1}\in P
\]
has Levi component \((1,h)\) and modulus \(1\). Thus the action is
\begin{equation}\label{eq:GL5-affine-action}
 (h\Phi)(x)=\delta(h)\Phi(h^{-1}x).
\end{equation}
Evaluation at \(0\) is the quotient map \(\mathcal B\to\delta|_H\) with kernel \(\mathcal A\).

If \(\delta\) is supercuspidal, Proposition~\ref{prp:even-supercuspidal-vanishing} and \eqref{eq:GL5-orbit-Hom} give \(\Hom_H(\mathcal A,\one)=0\), and \eqref{eq:GL5-restriction-injections} finishes the proof. Otherwise write
\[
 \delta=\St\bigl([\nu^{-(k-1)/2}\rho,\ldots,\nu^{(k-1)/2}\rho]\bigr),
 \qquad k\in\{2,4\},\qquad \deg\rho=4/k,
\]
with \(\rho\) unitary supercuspidal. The proof of Theorem~\ref{thm:intro-even-counterexample} for \(k=2\) and Remark~\ref{rmk:GL4-Steinberg-multiplicity} for \(k=4\), together with \eqref{eq:GL5-orbit-Hom} and \eqref{eq:GL5-restriction-injections}, give
\begin{equation}\label{eq:GL5-multiplicity-bound}
 \dim\Hom_H(\mathcal B,\one)
 \leq\dim\Hom_{P(\Sp_4)}(\delta,\one)\leq1.
\end{equation}

For \(a_t=\operatorname{diag}(t,1,1,t^{-1})\) and \(\ell\in\Hom_{P(\Sp_4)}(\delta,\one)\), we will prove
\begin{equation}\label{eq:GL5-period-torus}
 \ell(\delta(a_t)v)=|t|^{k+2}\ell(v)
 \qquad(t\in F^\times,\ v\in V_\delta).
\end{equation}

We first reduce the calculation to a compactly induced representation. Retain \(U,U',S,\xi,N\) from the proof of Theorem~\ref{thm:intro-even-counterexample}, and set \(E=\delta_{U'}\) and \(\tau=\delta_{U,\xi}\). Since \(P(\Sp_4)=S\ltimes U'\), the functional \(\ell\) descends to \(\bar\ell\in\Hom_S(E,\one)\). Applying \eqref{eq:rank-four-clifford-sequence} to \(\delta\) gives
\[
 0\longrightarrow M\simeq\cInd_N^S\tau\longrightarrow E\longrightarrow\delta_U\longrightarrow0.
\]
Equations~\eqref{eq:rank-four-Uprime}--\eqref{eq:rank-four-semi-Whittaker} for \(k=2\), and \eqref{eq:rank-four-semi-Whittaker 2} in Remark~\ref{rmk:GL4-Steinberg-multiplicity} for \(k=4\), give \(\Hom_S(\delta_U,\one)=0\) and \(\dim\Hom_N(\tau,\one)=1\). Thus restriction \(\Hom_S(E,\one)\to\Hom_S(M,\one)\) is injective. Since \(a_t\) normalizes \(U\), \(U'\), and \(P(\Sp_4)\), and commutes with \(S\), it acts on \(E\) and preserves \(M\). It therefore suffices to prove \(\bar\ell(a_tm)=|t|^{k+2}\bar\ell(m)\) for \(m\in M\).

We describe \(\bar\ell|_M\) using the model \(M\simeq\cInd_N^S\tau\). By \cite[\S3, proof of Proposition~1]{DijolsPrasad}, \(M\) consists of compactly supported sections on the nonzero character orbit of \(U/U'\), and the quotient map \(p_\xi:E\to\tau\) is evaluation at \(\xi\). The isomorphism in \cite[proof of Proposition~2.23]{Bern-Zele-Survey} is therefore
\begin{equation}\label{eq:GL5-Clifford-realization}
 m\longmapsto f_m,\qquad f_m(s)=p_\xi(sm),
\end{equation}
with \(f_m(ns)=\tau(n)f_m(s)\) and \(S\) acting by right translation.

Choose \(0\neq\lambda\in\Hom_N(\tau,\one)\). Since \(N\) and \(S\) are unimodular, Frobenius reciprocity \cite[proof of Proposition~2.29]{Bern-Zele-Survey} identifies \(\Hom_N(\tau,\one)\) with \(\Hom_S(M,\one)\) by
\[
 \lambda\longmapsto I_\lambda,\qquad
 I_\lambda(f)=\int_{N\backslash S}\lambda(f(s))\,d\dot s.
\]
Here \(d\dot s\) is induced by additive Haar measure on \(F^2\), identifying \(N\backslash S\) with \(F^2\setminus\{0\}\) by the bottom row of the middle \(2\times2\) block. The integrand descends to a locally constant, compactly supported function on \(N\backslash S\), so the integral converges absolutely. Since \(\Hom_N(\tau,\one)\) is one-dimensional, \(\bar\ell|_M\) is a multiple of \(m\mapsto I_\lambda(f_m)\). It therefore suffices to prove \(I_\lambda(f_{a_tm})=|t|^{k+2}I_\lambda(f_m)\).

For \(t\in F^\times\), set
\[
 s_t=\operatorname{diag}(1,t,t^{-1},1),\qquad
 d_t=a_ts_t=\operatorname{diag}(t,t,t^{-1},t^{-1}).
\]
The element \(d_t\) normalizes \(U\) and fixes \(\xi\), so it acts on \(\tau\). Formula~\eqref{eq:GL5-Clifford-realization} gives
\begin{equation}\label{eq:GL5-torus-on-M}
 f_{a_tm}(s)=p_\xi(sa_tm)
 =p_\xi(d_ts_t^{-1}sm)
 =\tau(d_t)f_m(s_t^{-1}s).
\end{equation}
Thus we must compute the effect of \(\tau(d_t)\) on \(\lambda\) and the change of measure under \(s\mapsto s_t^{-1}s\).

The pullback of \(\lambda\) is \((N_4,\psi_{2,2})\)-equivariant and therefore factors through \(\delta_{U_{2,2}}\), since \(\psi_{2,2}|_{U_{2,2}}=1\). Set \(D=\St([\nu^{-(k-2)/4}\rho,\ldots,\nu^{(k-2)/4}\rho])\). The Jacquet-module identities in the proof of Theorem~\ref{thm:intro-even-counterexample} and Remark~\ref{rmk:GL4-Steinberg-multiplicity} give
\begin{equation}\label{eq:GL5-Jacquet-22}
 \Jac_{P_{2,2}}^{\Gl_4}(\delta)
 \simeq\nu^{k/4}D\boxtimes\nu^{-k/4}D.
\end{equation}
On this normalized Jacquet module, the two central characters of \(D\) cancel at \(d_t\), while the twists contribute
\[
 |\det(tI_2)|^{k/4}|\det(t^{-1}I_2)|^{-k/4}=|t|^k.
\]
Passing to unnormalized coinvariants contributes \(\delta_{P_{2,2}}(d_t)^{1/2}=|t|^4\). Thus
\begin{equation}\label{eq:GL5-lambda-torus}
 \lambda(\tau(d_t)w)=|t|^{k+4}\lambda(w).
\end{equation}

For the change of measure, \(s_t\) normalizes \(N\), and \(s\mapsto s_t^{-1}s\) multiplies the bottom row by \(t\). Hence, for \(\varphi\in C_c^\infty(N\backslash S)\),
\[
 \int_{N\backslash S}\varphi(s_t^{-1}s)\,d\dot s
 =|t|^{-2}\int_{N\backslash S}\varphi(s)\,d\dot s.
\]
Combining this with \eqref{eq:GL5-torus-on-M} and \eqref{eq:GL5-lambda-torus}, we obtain
\[
 I_\lambda(f_{a_tm})=|t|^{k+4}|t|^{-2}I_\lambda(f_m)
 =|t|^{k+2}I_\lambda(f_m).
\]
Since \(\Hom_S(E,\one)\to\Hom_S(M,\one)\) is injective, the same
identity holds for \(\bar\ell\) on \(E\), proving \eqref{eq:GL5-period-torus}.

We now use Fourier transform on \(\mathcal B\) to prove that \(\Hom_H(\mathcal B,\one)=0\).
Suppose that \(0\neq L\in\Hom_H(\mathcal B,\one)\). By \eqref{eq:GL5-restriction-injections}, its restriction to \(\mathcal A\) is nonzero. The isomorphism for \(\mathcal A\) in \eqref{eq:GL5-orbit-filtration} is realized by \(\Phi\mapsto f_\Phi\), where \(f_\Phi(h)=\delta(h)\Phi(h^{-1}e_1)\). Since \(H\) and \(P(\Sp_4)\) are unimodular, Frobenius reciprocity \cite[proof of Proposition~2.29]{Bern-Zele-Survey} identifies \(L|_{\mathcal A}\) with a nonzero \(\ell\in\Hom_{P(\Sp_4)}(\delta,\one)\). The same proof expresses \(L|_{\mathcal A}\) as an integral over \(P(\Sp_4)\backslash H\), as follows.

Identify \(P(\Sp_4)\backslash H\) with \(W^\times=W\setminus\{0\}\) by \(P(\Sp_4)h\mapsto h^{-1}e_1\), and normalize the quotient measure \(d\dot h\) to correspond to the \(H\)-invariant additive Haar measure \(dx\) on \(W\). For each \(x\in W^\times\), choose \(h_x\in H\) satisfying \(h_xe_1=x\). Then
\begin{equation}\label{eq:GL5-open-functional}
 L|_{\mathcal A}(\Phi)
 =\int_{P(\Sp_4)\backslash H}\ell\bigl(f_\Phi(h)\bigr)\,d\dot h
 =\int_{W^\times}\ell\bigl(\delta(h_x^{-1})\Phi(x)\bigr)\,dx,
 \qquad \Phi\in\mathcal A.
\end{equation}
The expression is independent of the choice of \(h_x\) by the \(P(\Sp_4)\)-invariance of \(\ell\).

For \(t\in F^\times\), let \((D_t\Phi)(x)=\Phi(tx)\). This operator commutes with \(H\) and preserves both \(\mathcal A\) and \(\mathcal B\). In \eqref{eq:GL5-open-functional}, substitute \(y=tx\) and choose \(h_{t^{-1}y}=h_ya_t^{-1}\). Then \eqref{eq:GL5-period-torus} gives
\begin{align*}
 L(D_t\Phi)
 &=|t|^{-4}\int_{W\setminus\{0\}}
   \ell\bigl(\delta(a_t)\delta(h_y^{-1})\Phi(y)\bigr)\,dy\\
 &=|t|^{k-2}L(\Phi),\qquad \Phi\in\mathcal A.
\end{align*}
The \(H\)-invariant functional \(LD_t-|t|^{k-2}L\) on \(\mathcal B\) therefore vanishes on \(\mathcal A\). By \eqref{eq:GL5-orbit-Hom}, it vanishes on all of \(\mathcal B\). Hence
\begin{equation}\label{eq:GL5-dilation-on-B}
 LD_t=|t|^{k-2}L\qquad(t\in F^\times).
\end{equation}

Choose the Haar measure on \(W\) self-dual with respect to \(\psi\) and the pairing \(\langle x,y\rangle={}^txJ_4y\). Define
\begin{equation}\label{eq:GL5-Fourier}
 (\mathcal F\Phi)(x)=\int_W\Phi(y)\psi(\langle x,y\rangle)\,dy.
\end{equation}
This defines an automorphism of \(\mathcal B=C_c^\infty(W)\otimes_{\mathbb C}V_\delta\), acting by scalar Fourier transform on the first factor. With our alternating pairing and self-dual measure, Fourier inversion gives \(\mathcal F^2=1\); see \cite[\S4.5]{PoonenTate}. For \(h\in H\), we have \(\det h=1\) and \(\langle x,hz\rangle=\langle h^{-1}x,z\rangle\). Thus \eqref{eq:GL5-affine-action} and the substitution \(y=hz\) give
\begin{align*}
 (\mathcal F(h\Phi))(x)
 &=\delta(h)\int_W\Phi(h^{-1}y)\psi(\langle x,y\rangle)\,dy\\
 &=\delta(h)\int_W\Phi(z)\psi(\langle h^{-1}x,z\rangle)\,dz
 =(h\mathcal F\Phi)(x).
\end{align*}
Thus \(\mathcal F\) is \(H\)-equivariant and \(L\mathcal F\) is a nonzero \(H\)-invariant functional on \(\mathcal B\). By \eqref{eq:GL5-multiplicity-bound},
\begin{equation}\label{eq:GL5-Fourier-on-period}
 L\mathcal F=cL\qquad\text{for some }c\in\mathbb C^\times.
\end{equation}
Finally, the substitution \(z=ty\) in \eqref{eq:GL5-Fourier} gives
\[
 (\mathcal FD_t\Phi)(x)
 =|t|^{-4}(\mathcal F\Phi)(t^{-1}x),
 \qquad\text{that is,}\qquad
 \mathcal FD_t=|t|^{-4}D_{t^{-1}}\mathcal F.
\]
Combining this with \eqref{eq:GL5-dilation-on-B} and \eqref{eq:GL5-Fourier-on-period}, we obtain
\[
 c|t|^{k-2}L
 =L\mathcal FD_t
 =|t|^{-4}LD_{t^{-1}}\mathcal F
 =c|t|^{-k-2}L.
\]
Since \(cL\neq0\), this forces \(|t|^{2k}=1\) for every \(t\in F^\times\), a contradiction on taking \(t\) to be a uniformizer. Therefore \(\Hom_H(\mathcal B,\one)=0\), and \eqref{eq:GL5-restriction-injections} proves \eqref{eq:GL5-chi-delta-vanishing}. The initial reduction completes the proof of the necessary condition.
\end{proof}

Theorem~\ref{thm:intro-GL5-necessary} settles the necessary condition for
\(n=2\): the extra invariant functional supplied by
Theorem~\ref{thm:intro-even-counterexample} on the open orbit does not extend
to the induced representation. Together with the sufficient condition in
Theorem~\ref{thm:intro-sufficient}, this leads to the general statement
posed in the introduction.

\begin{restatedconjectureinner}
Let \(n\geq3\), and let \(\pi\) be an irreducible tempered representation
of \(\Gl_{2n+1}(F)\). If
\(\Hom_{\Sp_{2n}(F)}(\pi,\one)\neq0\), then
\[
 \pi\simeq\sigma_n\times\sigma_{n+1}
\]
for irreducible tempered representations \(\sigma_n\) of \(\Gl_n(F)\)
and \(\sigma_{n+1}\) of \(\Gl_{n+1}(F)\).
\end{restatedconjectureinner}

\par\vspace{0.5cm}
\appendix
\section{Parabolic orbits and stabilizers}
\label{app:parabolic-orbits}

We describe the parabolic orbits, their stabilizers and their images in the Levi subgroup for \(\Sp_{2n}\backslash\Gl_{2n+1}\), extending the calculations for \(P_{n,n+1}\) in \cite[\S5]{Xi-Lu} to \(P_{l,2n+1-l}\), \(1\leq l\leq n\).

Let \(G=\Gl_{2n+1}(F)\), let \(H=\Sp_{2n}(F)\) be embedded by \(h\mapsto\operatorname{diag}(h,1)\), and put \(P=P_{l,2n+1-l}\) and \(Q=Q_{2n+1}\). Write
\[
u(x)=\begin{pmatrix}I_{2n}&x\\&1\end{pmatrix},\qquad x\in F^{2n},
\]
and let \(\omega\) be the permutation matrix corresponding to \((l,2n+1)\). There are two \((P,Q)\)-double cosets, giving \(G=PQ\sqcup P\omega Q\), with
\[
P\cap Q=P_{l,2n-l}\ltimes u(F^{2n}),\qquad
Q\cap\omega^{-1}P\omega
 =P_{l-1,2n-l+1}\ltimes u(E_{l-1}),
\]
where \(E_k=\Span\{e_1,\ldots,e_k\}\). We use the conventions
\(E_0=0\), \(P_{0,2n}=\Gl_{2n}\), and \(\Sp_0=\{1\}\).
Empty spans are understood to be zero, and matrix blocks of size zero are omitted.

\subsection{Symplectic orbits on the Grassmannian}

For \(0\leq k\leq n\), the map \(P_{k,2n-k}g^{-1}\mapsto gE_k\) identifies \(P_{k,2n-k}\backslash\Gl_{2n}\) with \(\operatorname{Gr}(k,F^{2n})\). For a subspace \(W\), let \(W_d=W\cap W^\perp\) be its radical. By Witt's theorem, its \(\Sp_{2n}(F)\)-orbit is determined by \(2a=\dim(W/W_d)\). The possible values are \(0\leq a\leq\floor*{k/2}\), and we may take the representatives
\[
V_{k,a}
 =\Span\{e_1,\ldots,e_{k-a},f_a,\ldots,f_1\}.
\]
Let \(w_{k,a}\) be the permutation matrix corresponding to
\[
(k,2n)(k-1,2n-1)\cdots(k-a+1,2n-a+1),
\qquad w_{k,0}=1.
\]
Thus
\begin{equation}\label{eq:symplectic-grassmannian-orbits}
\Gl_{2n}
 =\bigsqcup_{a=0}^{\floor*{k/2}}
 P_{k,2n-k}w_{k,a}\Sp_{2n}(F).
\end{equation}

The stabilizer of \(V_{k,a}\) is \(S_{k,a}:=\Sp_{2n}(F)\cap w_{k,a}^{-1}P_{k,2n-k}w_{k,a}\). Put \(r=k-2a\) and \(m=2n-2k+2a\). Relative to the block sizes \((a,r,m,r,a)\), its elements have the form
\begin{equation}\label{eq:base-stabilizer}
\begin{pmatrix}
A&0&0&Y_1&B\\
Y_2&g&X&Z&M_2\\
0&0&h&K&0\\
0&0&0&g^\#&0\\
C&0&0&M_1&D
\end{pmatrix},
\qquad
\begin{pmatrix}A&B\\C&D\end{pmatrix}\in\Sp_{2a}(F),\quad
g\in\Gl_r(F),\quad h\in\Sp_m(F),
\end{equation}
where
\[
g^\#=\w_{0,r}\,{}^tg^{-1}\w_{0,r},\qquad
X^\natural=J_m\,{}^tX\w_{0,r},\qquad
K=hX^\natural g^\#.
\]
The blocks \(Y_2,M_2,X\) may be chosen freely. The symplectic relations then determine \(Y_1,M_1\) and the skew-symmetric part of \(\w_{0,r}g^{-1}Z\). More precisely, writing
\[
\Theta
 ={}^tY_1\w_{0,a}M_1
  -{}^tM_1\w_{0,a}Y_1
  +{}^tKJ_mK,
\]
we have
\[
{}^tg^\#\w_{0,r}Z-{}^tZ\w_{0,r}g^\#=\Theta,
\qquad
Z=g\w_{0,r}\left(T+\frac12\Theta\right),
\quad {}^tT=T.
\]
A Levi subgroup of \(S_{k,a}\) is
\begin{equation}\label{eq:base-stabilizer-levi}
\Sp_{2a}(F)\times\Gl_r(F)\times\Sp_m(F).
\end{equation}
Let
\[
Y_{k,a}
 =\Span\{e_{k-a+1},\ldots,e_n,
          f_n,\ldots,f_{k-a+1}\},
\quad
R_{k,a}^\vee
 =\Span\{f_{k-a},\ldots,f_{a+1}\},
\]
of dimensions \(m\) and \(r\), respectively. We identify \(F^{2n}/V_{k,a}\) with \(Y_{k,a}\oplus R_{k,a}^\vee\); the induced action of \(S_{k,a}\) is
\begin{equation}\label{eq:base-stabilizer-quotient-action}
(y,\eta)\longmapsto(hy+K\eta,g^\#\eta).
\end{equation}

For the Levi projections below, we use the bases \((e_1,\ldots,e_a,f_a,\ldots,f_1)\), \((e_{a+1},\ldots,e_{k-a})\) and \((f_{k-a},\ldots,f_{a+1})\) for the \(\Sp_{2a}\), \(g\) and \(g^\#\) blocks, respectively, and the displayed basis of \(Y_{k,a}\) for the \(h\) block. These bases are transported by the orbit representative and reordered as in each formula. The additional one-dimensional block corresponds to \(e_{2n+1}\) for orbits in \(PQ\), and to \(e_{2n+1}-v_{i,a}\) for orbits in \(P\omega Q\). In \eqref{eq:M3a}, we also reverse both \(r\)-dimensional bases so that the first block belongs to \(Q_r\).

\subsection{Orbits in \texorpdfstring{\(PQ\)}{PQ}}

For \(0\leq a\leq\floor*{l/2}\), put \(g_a^1=\operatorname{diag}(w_{l,a},1)\) and \(\cO_a^1=Pg_a^1H\). Applying \eqref{eq:symplectic-grassmannian-orbits} with \(k=l\) gives
\begin{equation}\label{eq:first-cell-orbits}
PQ=\bigsqcup_{a=0}^{\floor*{l/2}}\cO_a^1.
\end{equation}
The stabilizers in \(H\) and \(P\) are \(H_a^1=\operatorname{diag}(S_{l,a},1)\) and \(P_a^1=g_a^1H_a^1(g_a^1)^{-1}\), respectively. Let \(M\simeq\Gl_l(F)\times\Gl_{2n+1-l}(F)\) be the standard Levi subgroup of \(P\). With \(r=l-2a\) and \(m=2n-2l+2a\), the image \(M_a^1=\operatorname{pr}_M(P_a^1)\) is conjugate in \(M\) to
\begin{equation}\label{eq:first-cell-levi-projection}
\left\{
\left(
\begin{pmatrix}s&0\\U&g\end{pmatrix},
\begin{pmatrix}h&K&0\\0&g^\#&0\\0&0&1\end{pmatrix}
\right):
\begin{array}{c}
s\in\Sp_{2a}(F),\ g\in\Gl_r(F),\ h\in\Sp_m(F),\\
U\in\Mat_{r,2a}(F),\ K\in\Mat_{m,r}(F)
\end{array}
\right\}.
\end{equation}
For the identity orbit with \(l=n\), this gives
\[
\operatorname{pr}_M(P\cap H)
 =\{(g,\operatorname{diag}(g^\#,1)):g\in\Gl_n(F)\},
\]
as in Proposition~\ref{prp:odd-closed-orbit-period}.

\subsection{Orbits in \texorpdfstring{\(P\omega Q\)}{P omega Q}}

Put \(k=l-1\), \(\widetilde w_a=\operatorname{diag}(w_{k,a},1)\), and \(L=Q\cap\omega^{-1}P\omega=P_{k,2n-k}\ltimes u(E_k)\). The natural map
\[
L\backslash Q/H\longrightarrow
P_{k,2n-k}\backslash\Gl_{2n}/H
\]
has fiber \(S_{k,a}\backslash(F^{2n}/V_{k,a})\) over the class of \(w_{k,a}\). We may take the representatives
\[
v_{1,a}=0,\qquad
v_{2,a}=e_{k-a+1},\qquad
v_{3,a}=f_{a+1},
\]
where the third occurs only when \(r=k-2a>0\). These exhaust the orbits by \eqref{eq:base-stabilizer-quotient-action} and the transitivity of \(\Sp_m\) and \(\Gl_r\) on nonzero vectors.

Set
\[
x_{1,a}=0,\qquad
x_{2,a}=w_{k,a}v_{2,a}
 =\begin{cases}
   e_{k+1},&a=0,\\
   f_a,&a>0,
  \end{cases}
\qquad
x_{3,a}=f_{a+1},
\]
and put \(g_{i,a}=\omega u(x_{i,a})\widetilde w_a\) and \(\cO_{i,a}=Pg_{i,a}H\). Together with the orbits in \(PQ\), these give
\begin{equation}\label{eq:odd-orbit-decomposition}
\begin{aligned}
P\backslash G/H={}&
 \{\cO_a^1:0\leq a\leq\floor*{l/2}\}\\
&\sqcup
 \{\cO_{i,a}:i=1,2,\
       0\leq a\leq\floor*{(l-1)/2}\}\\
&\sqcup
 \{\cO_{3,a}:0\leq a\leq\floor*{(l-2)/2}\}.
\end{aligned}
\end{equation}

The stabilizers in \(H\) are
\[
H_{i,a}
 =\{s\in S_{k,a}:sv_{i,a}-v_{i,a}\in V_{k,a}\}.
\]
Thus \(H_{1,a}=S_{k,a}\). In the coordinates of \eqref{eq:base-stabilizer}, the condition for \(H_{2,a}\) is \(hv_{2,a}=v_{2,a}\), so the \(h\) block lies in a conjugate of \(P(\Sp_m)\). For \(r>0\),
\[
H_{3,a}
 =\{s:g^\#v_{3,a}=v_{3,a},\ Kv_{3,a}=0\}.
\]
Writing \(R_{k,a}=Fe_{a+1}\oplus\Span\{e_{a+2},\ldots,e_{k-a}\}\), these conditions become
\begin{equation}\label{eq:third-orbit-stabilizer-condition}
g=\begin{pmatrix}1&0\\ \gamma&g_{r-1}\end{pmatrix},
\qquad
X=\begin{pmatrix}0_{1,m}\\X_{r-1,m}\end{pmatrix}.
\end{equation}
The stabilizer in \(P\) is
\begin{equation}\label{eq:odd-orbit-stabilizers}
\begin{aligned}
P_{i,a}
 &=g_{i,a}H_{i,a}g_{i,a}^{-1}\\
 &=\omega\left\{
 \begin{pmatrix}
 w_{k,a}s w_{k,a}^{-1}&
 (I_{2n}-w_{k,a}s w_{k,a}^{-1})x_{i,a}\\
 0&1
 \end{pmatrix}:s\in H_{i,a}
 \right\}\omega^{-1}.
\end{aligned}
\end{equation}

\subsection{Levi projections}

Let \(M_{i,a}=\operatorname{pr}_M(P_{i,a})\), and write \(r=l-2a-1\) and \(m=2n-2l+2a+2\). Projecting \eqref{eq:odd-orbit-stabilizers} to \(M\) and using the symplectic relations gives the following descriptions, up to conjugacy in \(M\). For \(0\leq a\leq\floor*{(l-1)/2}\), we obtain
\begin{align}
M_{1,a}
&=
\left\{
\left(
\begin{pmatrix}s&0&0\\U&g&0\\0&0&1\end{pmatrix},
\begin{pmatrix}h&K\\0&g^\#\end{pmatrix}
\right):
\begin{array}{c}
s\in\Sp_{2a}(F),\ g\in\Gl_r(F),\ h\in\Sp_m(F),\\
U\in\Mat_{r,2a}(F),\ K\in\Mat_{m,r}(F)
\end{array}
\right\},
\label{eq:M1a}\\[3pt]
M_{2,a}
&=
\left\{
\left(
\begin{pmatrix}
g&U&g\w_{0,r}\,{}^tK^-\\
0&s&0\\
0&0&1
\end{pmatrix},
\begin{pmatrix}h_+&K\\0&g^\#\end{pmatrix}
\right):
\begin{array}{c}
s\in\Sp_{2a}(F),\ g\in\Gl_r(F),\
h_+\in P(\Sp_m),\\
U\in\Mat_{r,2a}(F),\ K\in\Mat_{m,r}(F)
\end{array}
\right\}.
\label{eq:M2a}
\end{align}
Here \(K^-\in\Mat_{1,r}(F)\) denotes the last row of \(K\) in the symplectic basis for which \(h_+\) fixes the first vector. The upper-right block in the first matrix of \eqref{eq:M2a} comes from the relation \(-X_{\bullet,1}=g\w_{0,r}{}^tK^-\), where \(X_{\bullet,1}\) is the first column of \(X\) in \eqref{eq:base-stabilizer}, before reordering.

Finally, for \(0\leq a\leq\floor*{(l-2)/2}\), so that \(r>0\), we have
\begin{equation}\label{eq:M3a}
M_{3,a}
=
\left\{
\left(
\begin{pmatrix}
q&U&Z\\
0&s&\lambda\\
0&0&1
\end{pmatrix},
\begin{pmatrix}
h&\begin{matrix}0_{m,1}&K\end{matrix}\\
0&q^\#
\end{pmatrix}
\right):
\begin{array}{c}
q\in Q_r,\ s\in\Sp_{2a}(F),\ h\in\Sp_m(F),\\
U\in\Mat_{r,2a}(F),\ Z\in\Mat_{r,1}(F),\\
\lambda\in\Mat_{2a,1}(F),\
K\in\Mat_{m,r-1}(F),\\
{}^tsJ_{2a}\lambda={}^tUe_r
\end{array}
\right\},
\end{equation}
where \(q^\#=\w_{0,r}{}^tq^{-1}\w_{0,r}\). The formulas \eqref{eq:M1a} and \eqref{eq:M2a} remain valid when \(r=0\), with the empty blocks omitted; in \eqref{eq:M3a}, we use \(Q_1=Q_1^\#=\{1\}\) when \(r=1\).

\bibliographystyle{alpha}
\bibliography{refs}
\end{document}